\documentclass[11pt]{article}

\usepackage{amsfonts} 
\usepackage{amsbsy} 

\newcommand{\ih}{\'{\i}}
\newcommand{\eh}{\hspace{.06in}}
\newcommand{\C}{\mathbb{C}} 
\newcommand{\E}{\mathbb{E}} 
\newcommand{\N}{\mathbb{N}} 
\newcommand{\R}{\mathbb{R}} 
\newcommand{\Z}{\mathbb{Z}} 
\newcommand{\m}{\frac{1}{2}}   
\newcommand{\see}{\eh{\Leftrightarrow}\eh}  
\newcommand{\ds}{\displaystyle}  
\newcommand{\BE}{\begin{equation}}  
\newcommand{\EE}{\end{equation}}
\newcommand{\tg}{\tilde{g}}
\newcommand{\tp}{\tilde{p}}
\newcommand{\tq}{\tilde{q}}
\newcommand{\tr}{\tilde{r}}
\newcommand{\tJ}{\tilde{J}}  
\newcommand{\tM}{\tilde{M}}
\newcommand{\tX}{\tilde{X}}
\newcommand{\af}{\alpha}   
\newcommand{\bs}{\bf s}
\newcommand{\brho}{\boldsymbol{\rho}}
\newcommand{\btau}{\boldsymbol{\tau}}
\newcommand{\eps}{\varepsilon}
\newcommand{\A}{\mathcal A}
\newcommand{\B}{\mathcal B}
\newcommand{\Cc}{\mathcal C}
\newcommand{\D}{\mathcal D}
\newcommand{\Dd}{\mathsf D}
\newcommand{\Ee}{\mathcal E}
\newcommand{\F}{\mathcal F}
\newcommand{\G}{\mathcal G}
\newcommand{\Nn}{\mathcal N} 
\newcommand{\Ns}{\mathsf N} 
\newcommand{\Rr}{\mathcal R}
\newcommand{\Ss}{\mathcal S}
\newcommand{\T}{\mathcal T}
\newcommand{\deh}{\partial}
\newcommand{\fns}{\footnotesize} 
\newcommand{\Lim}[1]{\lower5pt\hbox{${{\ds\lim}\atop #1}$}}

\begin{document}
\centerline{\large{\bf A CHARACTERISATION OF THE}}
\centerline{\large{\bf HOFFMAN-WOHLGEMUTH SURFACES IN}}
\centerline{\large{\bf TERMS OF THEIR SYMMETRIES}}
\bigskip

\centerline{P{\fns INIO} S{\fns IM\~OES} \& V{\fns AL\'ERIO} R{\fns AMOS} B{\fns ATISTA}}
\bigskip

\begin{abstract}
For an embedded singly periodic minimal surface $\tilde{M}$ with genus $\varrho\ge4$ and annular ends, some weak symmetry hypotheses imply its congruence with one of the Hoffman-Wohlgemuth examples. We give a very geometrical proof of this fact, along which they come out many valuable clues for the understanding of these surfaces.
\end{abstract}
\ \\
\\
\centerline{\bf 1. Introduction}
\

The beauty of a characterisation theorem resides particularly in its demonstration, where a lot of intrinsic and fascinating properties are revealed. For complete embedded minimal surfaces of finite total curvature in the euclidean space $\E=\R^3$, R.Schoen published a strong result in 1983: if $S$ is such a surface, then it must be the catenoid providing it has exactly two ends (see [{\bf 28}]). Later in 1991, F.L\'opez and A.Ros proved that, if $S$ has genus zero, then it is the catenoid or the flat plane (see [{\bf 14}]). Together, their works showed that other examples of $S$ should have positive genera {\it and} more than two ends. Meanwhile, C.Costa characterised all minimal tori $S$ with three ends (see [{\bf 1}]), but a torus $S$ with four ends or more could not exist by the Hoffman-Meeks' conjecture that $\sharp$ ends $\le$ genus+2. Higher genus examples can be found in [{\bf 31}].
\\

In 1990, D.Hoffman and W.Meeks gave examples of $S$ with three ends and arbitrary positive genus, which in 1995 were generalised by D.Hoffman and H.Karcher (see [{\bf 9}] and [{\bf 7}]). Under symmetry hypothesis, in 2001 F.Mart\ih n and M.Weber classified them (see [{\bf 17}]). One year later, M.Traizet replaced the symmetry hypothesis by the weaker concept of {\it configuration} and got a characterisation of the Hoffman-Karcher two-parameter family (see [{\bf 29}]). Moreover, in the same work he gave examples of totally asymmetric $S$, answering the open question from [{\bf 7}, sec 5.2]. 
\\

Traizet's surfaces have high fixed genus, 5 ends and can assume different configurations. They show that any classification result of $S$ will need more constraints, for instance, a fixed conformal structure. However, in [{\bf 27}] it is shown that, when self-intersections are allowed at the ends of $S$, the conformal structure, even together with symmetry constraint, is insufficient to characterise the surface.
\\

After having discussed $S$ in $\E=\R^3$, it is important to mention the advances in $\E=\R^3/\T$, where $\T$ is a cyclic translation group. If $S$ is a torus with a finite number of planar ends, then $S$ belongs to Riemann's family according to [{\bf 19}]. The genus-one hypothesis is necessary because of an unpublished work from F.Wei, {\it Adding handles to the Riemann examples}. However, further characterisation results had to impose more constraints to be accomplished. In 1997 and 2000, the beautiful works [{\bf 15}], [{\bf 16}] from F.Mart\ih n and D.Rodr\ih guez showed that mild hypotheses on ends, genus and symmetries imply that $S$ is one of the Callahan-Hoffman-Meeks' examples [{\bf 3}]. The symmetry conditions are necessary due to a work from M.Callahan, D.Hoffman and H.Karcher (see [{\bf 2}]). For results on $S$ with helicoidal and Scherk-ends, see [{\bf 20}] and [{\bf 25}]. 
\\
\input epsf
\begin{figure}[ht]
\centerline{
\epsfxsize 7.5cm
\epsfbox{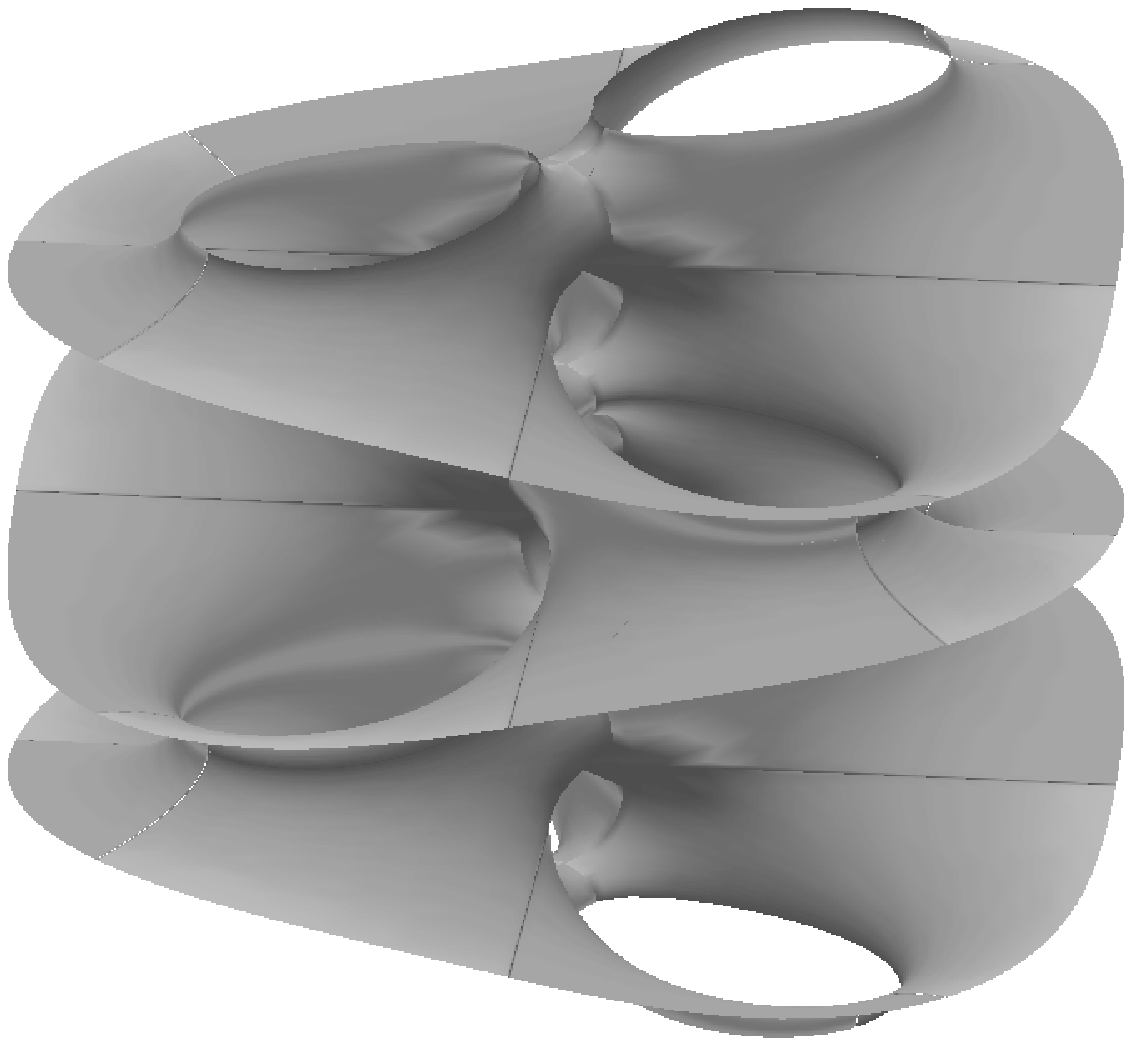}
\epsfxsize 8cm
\epsfbox{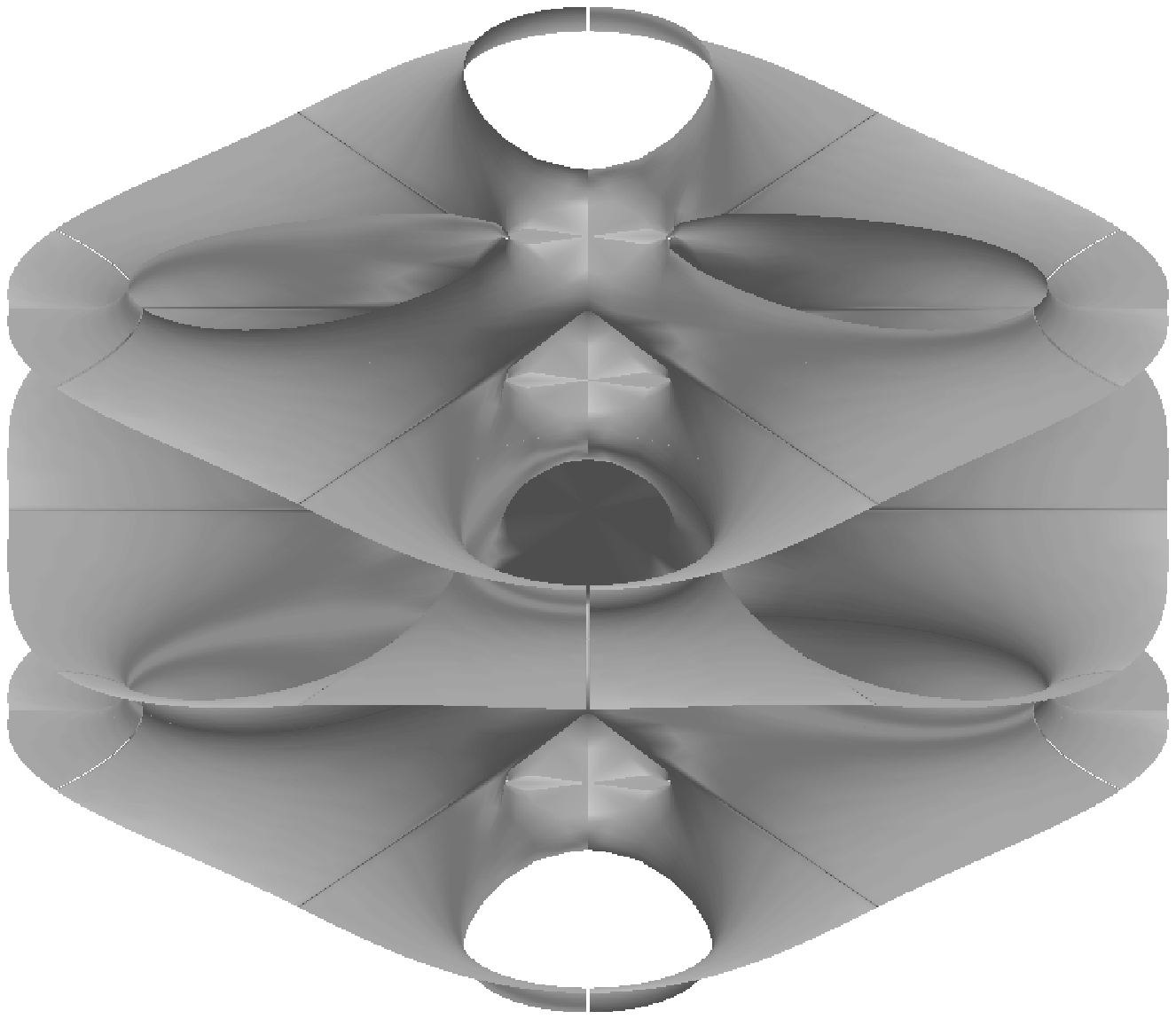}}
\caption{The Hoffman-Wohlgemuth surfaces of genera 5 and 9.}
\end{figure}

If $S$ is doubly periodic, the reader will find beautiful works like [{\bf 12}] and [{\bf 24}]. Back to the singly periodic case for $S$ with annular ends, our attention will now focus on the unpublished work from D.Hoffman and M.Wohlgemuth, {\it New embedded minimal surfaces of Riemann type}. These surfaces were obtained by adding Noevius handles to the examples in [{\bf 3}]. Of course, the sole addition of handles is by itself of little interest nowadays, except in the case of general results as [{\bf 11}]. Therefore, a characterisation theorem brings much more of new and good knowledge, particularly in the case of the Hoffman-Wohlgemuth surfaces. This present work is strongly inspired in the beautiful ideas of [{\bf 16}], but there are substantial differences, mainly because they deal with one-dimensional period problems, whereas the periods are two-dimensional in our case. In their work, the first part uses genus, ends and symmetry hypotheses to get Weierstrass data, and these allow 3 different family of surfaces. In the second part, hard computations of elliptic integrals finally show that just one family admits an embedded member, and only one, for any fixed odd genus starting from 3.
\\

Our first part is similar to theirs, but from the Weierstrass data $(g,dh)$ one gets 32 different families. However, simple geometric arguments quickly drop this number to 4. In the second part, a very basic handling of $\int gdh$ and $\int dh/g$ shows that, on 3 of the cases, the period is always open on a suitable closed curve. This is quite unexpected for two-dimensional problems, where non-existence in general follows from periods that can be separately solved, and then it lacks a {\it simultaneous} solution. Moreover, in our cases neither $\int gdh$ nor $\int dh/g$ will need any {\it explicit} formulation.
\\

The fact that one of the periods never closes is apparently due to the presence of a ``Gaussian geodesic''. By this concept we mean a planar curve of reflectional symmetry, which is the graph of an even real-analytic function $f:\R\to(0,1]$, where $f(0)=1$, $f'\ne 0$ in $\R^*$ and $\Lim{x\to\infty}{f(x)}=0$. Since 1997, when the second author started his doctoral studies in Germany, he observed that they failed all construction attempts of minimal surfaces containing a Gaussian geodesic. In total one tried 15 different examples and periods never closed. The same held for ``inverted Gaussians'', now with odd $f:\R\to(-1,1)$, $f'=0$ only at $0$ and $\Lim{x\to+\infty}{f(x)}=1$.
\\
\input epsf
\begin{figure}[ht]
\centerline{
\epsfxsize 15cm
\epsfbox{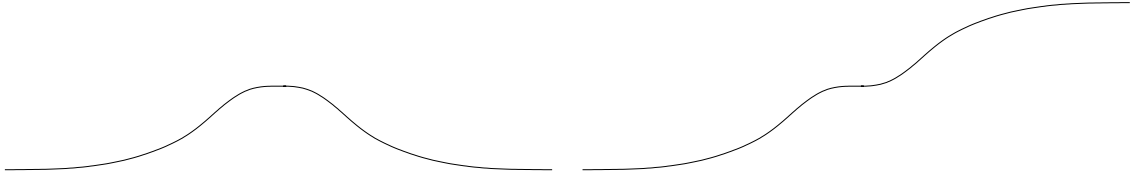}}
\caption{A standard and an inverted Gaussian.}
\end{figure}

Hitherto, it remains open the question whether an $S$ with Gaussian geodesic exists. However, Section 6 of this present work might bring some good ideas for a future study of this question. At last, the closure of periods succeeds for the Hoffman-Wohlgemuth surfaces, and yet the proof {\it is} easy (see Section 7). However, a {\it unique} solution could only be verified with numerical computation. This is typical for two-dimensional problems involving no L\'opez-Ros parameter, for till nowadays there is {\it just one} formal uniqueness demonstration of this kind, recently achieved by L.Ferrer and F.Mart\ih n (see [{\bf 6}]). Moreover, the proof in [{\bf 6}] is very laborious and reports a surface found 6 years beforehand (see [{\bf 8}]). Our result must be then interpreted in the sense that the Hoffman-Wohlgemuth family might occasionally admit two or more members with the same genus. 
\\

We shall give now a definition concerning {\it ends} of a complete Riemannian surface:
\\

{\bf Definition 1.1.} Consider a complete Riemannian surface $R$ and a sequence of enclosed compact balls $B_n\subset R$, $n\in\N$, with $\cup_{n=0}^\infty B_n=R$. Suppose there is $n_0>0$ and a connected component $\Ee$ of $R\setminus B_{n_0}$ such that $\Ee\setminus B_m$ and $\Ee\setminus B_n$ are homeomorphic for any $m$, $n\ge n_0$. In this case one has an equivalence relation $\Ee\setminus B_m\sim\Ee\setminus B_n$ and the corresponding equivalence class is called an \it end \rm of $S$. We also denote any of this class representative by $\Ee$. 
\\

Now we present the main theorem of this paper:
\\

{\bf Theorem 1.1.} \it Let $\tM$ be a properly embedded minimal surface in $\R^3$ such that 

i) $\tM$ has an infinite number of annular ends;

ii) $\tM$ is invariant under a cyclic group of screw motions $\T_\theta$;

iii) $\tM/\T_\theta$ has genus $\varrho\ge 4$ and exactly two ends;

iv) $|Iso(\tM/\T_\theta)|\ge 2(\varrho+3)$;

v) $\exists$ a conformal generator $\sigma\in$ Iso$(\tM/\T_\theta)$, and $r_1,r_2\in\tM$ such that $\sigma(r_i)\ne r_i\notin$

\ \ \ $\{r_j,\sigma(r_j)\}=\sigma(\{r_j,\sigma(r_j)\})$, for $\{i,j\}=\{1,2\}$.
\\
Then $(\varrho+1)/2$ is odd and $\tM$ is the Hoffman-Wohlgemuth surface $M_\varrho$ of genus $\varrho$.\rm
\\
\\
REMARKS: $M_\varrho$ exists only for $\varrho=4k+1$, $k\in\N^*$, as we shall see in Sections 3 and 7. Although apparently excessive, hypothesis {\it (v)} is necessary. Indeed, all surfaces from Callahan-Hoffman-Meeks of genus $\varrho>3$ verify {\it (i)-(iv)}, but not {\it (v)}. It could be replaced by ``Iso$(\tM/\T)<4(\varrho+1)$'', but in the praxis upper bounds for isometry groups are hard to compute. By a ``screw motion'' we mean a rotation about an axis followed by a translation {\it not necessarily} in the axis direction. 
\\

The present work was supported by FAPESP grant number 05/00026-3.
\\
\\
\centerline{\bf 2. Preliminaries}
\

In this section we state some basic definitions and theorems. Throughout this work, surfaces are considered connected and regular. Details can be found in [{\bf 10}], [{\bf 13}], [{\bf 21}] and [{\bf 23}]. 
\\

{\bf Theorem 2.1.} \it Let $X:R\to\E$ be a complete isometric immersion of a Riemannian surface $R$ into a three-dimensional complete flat space $\E$. If $X$ is minimal and the total Gaussian curvature $\int_R K dA$ is finite, then $R$ is biholomorphic to a compact Riemann surface $\overline{R}$ punched at a finite number of points.\rm 
\\

{\bf Theorem 2.2.} (Weierstrass representation). \it Let $R$ be a Riemann surface, $g$ and $dh$ meromorphic function and 1-differential form on $R$, such that the zeros of $dh$ coincide with the poles and zeros of $g$. Suppose that $X:R\to\E$, given by
\BE
   X(p):=Re\int^p(\phi_1,\phi_2,\phi_3),\eh\eh where\eh\eh
   (\phi_1,\phi_2,\phi_3):=\m(g^{-1}-g,ig^{-1}+ig,2)dh,\label{W_rep}
\EE
is well-defined. Then $X$ is a conformal minimal immersion. Conversely, every conformal minimal immersion $X:R\to\E$ can be expressed as (1) for some meromorphic function $g$ and 1-form $dh$.\rm
\\

{\bf Definition 2.1.} The pair $(g,dh)$ is the \it Weierstrass data \rm and $\phi_1$, $\phi_2$, $\phi_3$ are the \it Weierstrass forms \rm on $R$ of the minimal immersion $X:R\to X(R)\subset\E$.
\\

{\bf Theorem 2.3.} \it Under the hypotheses of Theorems 2.1 and 2.2, the Weierstrass data $(g,dh)$ extend meromorphically on $\overline{R}$.\rm  
\\

{\bf Theorem 2.4.} (Callahan-Hoffman-Meeks [{\bf 4}]). \it Suppose $X:R\to\R^3$ is a proper minimal embedding with more than one end. If $X(R)$ has an infinite group of symmetries, then it is either a catenoid or has the following properties:

1. $\sharp \ ends(X(R))=\infty$;

2. there is a screw motion $\T_\theta$ in $\R^3$ such that $\T_\theta(X(R))=X(R)$;

3. all annular ends of $X(R)$ are flat; 

4. the total curvature of $S:=X(R)/\T_\theta$ is finite if, and only if $\pi_1(S)$ is finitely generated. In this case $\int_SKdA=2\pi[\chi(S)-\sharp ends(S)]$.\rm
\\

The function $g$ is the stereographic projection of the Gau\ss \ map $N:R\to S^2$ of the minimal immersion $X$. It is a covering map of $\hat\C$ and $\int_SKdA=-4\pi$deg$(g)$. These facts will be largely used throughout this work.
\\
\\
\\
\centerline{\bf 3. The Weierstrass data of $\tM$}
\    

Considering the hypotheses {\it (i)} and {\it (iii)} of Theorem 1.1, at least one end of $\tM/\T_\theta$ must be an annulus in $\R^3$. Since $\tM$ is proper, the other end must be unbounded, and also an annulus because $\tM$ is embedded. Therefore, all ends of $\tM$ are annular and hence flat by Theorem 2.4. Now we apply {\it (iii)} to Theorem 2.4 and conclude that $\tM/\T_\theta$ has total curvature $-4\pi(\varrho+1)$. From Theorem 2.1 it follows that $\tM/\T_\theta$ is biholomorphic to a compact Riemann surface $\overline{M}$ punched at two points, because of {\it (iii)}. We call them $p_1$ and $p_2$.
\\

Now define $M:=\overline{M}\setminus\{p_1,p_2\}$. From the converse of Theorem 2.2, we have on $M$ a Weierstrass pair $(g,dh)$ which extends meromorphically on $\overline{M}$ by Theorem 2.3. Notice that deg$(g)=\varrho+1$. Up to a rigid motion in $\R^3$, $g(p_1)=0$ and so $g(p_2)=\infty$ because of {\it Alexander duality}. Let us write $\T_\theta=\rho\circ\tau$, where $\rho$ is a rotation about $Ox_3$ and $\tau$ a translation. The same arguments from [{\bf 15}, p187-8] easily generalise for non-vertical $\tau$, and they imply that $\tM$ is invariant under $\btau:=\tau^2$, whence also invariant under $\brho:=\rho^2$. If $\bs=\brho\circ\btau$ and $m=$ ord$(\brho)$, then both $\tM/\!<\btau>$ and $\tM/\!<\bs>$ will have the same total curvature, since each of them is an $m$-sheeted branched covering of $\tM/\!<\rho,\btau>$.
\\

ASSERTION 1: $|Iso(\tM/\!<\btau>)|\ge|Iso(\tM/\!<\bs>)|$.
\\

{\it Proof}. We have $\btau=\brho^{m-1}\circ\bs$. Hence the map $<\brho>\eh\to\eh<\brho,\bs>\!/\!<\btau>$, given by $\brho^i\mapsto\brho^i\!<\btau>$, is an isomorphism. Therefore, $|\!<\brho,\bs>\!/\!<\btau>\!|=$ $m$. If $G=Iso(\tM)$, then 
\[
   \biggl|\frac{G}{<\brho,\bs>}\biggl|\cdot\biggl|\frac{<\brho,\bs>}{<\btau>}\biggl|=
   \biggl|\frac{G}{<\btau>}\biggl|.
\]
The map $<\brho>\eh\to\eh<\brho,\bs>\!/\!<\bs>$, given by $\brho^i\mapsto\brho^i\!<\bs>$, is an epimorphism. Therefore, $|\!<\brho,\bs>\!/\!<\bs>\!|\le$ $m$. Now
\[
   \biggl|\frac{G}{<\brho,\bs>}\biggl|\cdot\biggl|\frac{<\brho,\bs>}{<\bs>}\biggl|=
   \biggl|\frac{G}{<\bs>}\biggl|,
\]
whence $|\,G/\!<\bs>\!|\le$ $|\,G/\!<\btau>\!|$.

\hfill q.e.d.
\ \\

Precisely CLAIM 3 of [{\bf 15}, p189] {\it implies} that $\tau$ is vertical. By following the same ideas as in [{\bf 15}, p189-0], one sees that $\tM$ is invariant under $\T_\theta^{-1}\circ\tau\circ\T_\theta=\tau$. Now the previous arguments apply for $\rho$, $\tau$ and $\T_\theta$ in the place of $\brho$, $\btau$ and $\bs$, respectively. Therefore, we can rewrite Theorem 1.1 as follows: 
\\

{\bf Theorem 3.1.} \it Let $\tM$ be a properly embedded minimal surface in $\R^3$ such that 

{\sl i}) $\tM$ has an infinite number of annular ends;

{\sl ii}) $\tM$ is invariant under a cyclic group of vertical translations $\T=<\tau>$;

{\sl iii}) $\tM/\T$ has genus $\varrho\ge 4$ and exactly two ends;

{\sl iv}) $|Iso(\tM/\T)|\ge 2(\varrho+3)$;

{\sl v}) $\exists$ a conformal generator $\sigma\in$ Iso$(\tM/\T)$, and $r_1,r_2\in\tM$ such that $\sigma(r_i)\ne r_i\notin$

\ \ \ $\{r_j,\sigma(r_j)\}=\sigma(\{r_j,\sigma(r_j)\})$, for $\{i,j\}=\{1,2\}$.
\\
Then $(\varrho+1)/2$ is odd and $\tM$ is the Hoffman-Wohlgemuth surface $M_\varrho$ of genus $\varrho$.\rm
\\

At this point, we re-define $M:=\tM/\T$, $\overline{M}:=M\cup\{p_1,p_2\}$ and $(g,dh)$ on $\overline{M}$ given by Theorem 2.3. Now the same arguments from [{\bf 16}, p448] firstly imply that the group $\Delta$ of automorphisms of $\overline{M}$ has a cyclic subgroup $\G:=\{A\in\Delta:A\eh{\rm is\eh holomorphic\eh and}\eh A(p_1)=p_1\}$. Secondly, if $J$ is a generator of $\G$, there is a corresponding symmetry $\tJ$ of $\tM$ which fixes a point in space. From {\it ({\sl v})} and the fact that $\G$ is cyclic, we may take $J=\sigma$.
\\

By Hurwitz's theorem $\Delta$ is finite, and so is $\G$. Therefore, ord$(\tJ)$ is finite and equals ord$(J)=n$. Without loss of generality we consider $\tJ(0,0,0)=(0,0,0)$. The rigidity of $\tM$ (see [{\bf 5}]) and the fact that $\tJ$ has a discrete fixed-point set on $\tM$ (possibly empty) imply that $\tJ$ keeps fixed the vertical $x_3$-axis. Since $[\Delta:\G]\le 4$, from {\it ({\sl iv})} we have $n\ge(\varrho+3)/2$, $\varrho\ge 4$. From [{\bf 15}, p189], $\tJ$ is a $2\pi/n$-rotation around $Ox_3$ composed with a reflection in $Ox_1x_2$. Up to a homothety, $\tau(x)=x+(0,0,2)$.    
\\

As in [{\bf 16}, p449], one defines for $q\in\overline{M}$ the stabiliser $S_q=\{f\in\G:f(q)=q\}$ and the orbit $O_q=\{q,J(q),\dots,J^{n-1}(q)\}$. Since $n=\sharp O_q\cdot\sharp S_q$, the Riemann-Hurwitz formula for the branched covering $\zeta:\overline{M}\to\overline{M}/J$ gives
\BE
   n\cdot\chi(\overline{M}/J)=\chi(\overline{M})+\biggl(2n-2+\sum_{q\in M}(\sharp S_q-1)\biggl).
\EE 

In (2) the term $2n-2$ corresponds to $p_{1,2}$. For each $q$ with $\sharp S_q>1$, consider the set $O_q=O_{J(q)}=\dots=O_{J^{n-1}(q)}$. There are exactly $s$ disjoint sets like that, $s\in\N^*$, and for each set we call its cardinality $m_i$, $i=1,\dots,s$. One rewrites (2) as follows:
\BE
   n\cdot\chi(\overline{M}/J)=2n-2\varrho+\sum_{i=1}^s(n-m_i).
\EE 

Up to re-indexing, from {\it ({\sl v})} we have $m_1=m_2=2$ and $s\ge 2$. Since $2n-\varrho\ge 3$, one guarantees that $\chi(\overline{M}/J)\ge 2/n>0$. Therefore, $\overline{M}/J$ is the Riemann sphere with Euler characteristic $2$. Hence (3) simplifies to 
\BE
   \sum_{i=1}^s(n-m_i)=2\varrho.
\EE 

Since $\tJ^2$ is a rotation around $Ox_3$, we conclude that $m_i\le 2$, $\forall$ $i$. For a minimal surface invariant under a rotation $\Rr$ about $Ox_3$, in [{\bf 3}] one proves that any flat horizontal end will have an order of pole (or zero) for $g$ given by $j\cdot$ord$(\Rr)+1$, where $j$ is a certain positive integer. Moreover, any regular point with vertical normal will have an order of pole (or zero) for $g$ given by $j\cdot$ord$(\Rr)-1$. Hence ord$_{p_i}(g)\ge n/2+1$ and ord$_{r_i}(g)\ge n/2-1$. For $s=2$, it follows from (4) that $n=2+\varrho$. Hence deg$(g)\ge 1+\varrho/2+1+2(1+\varrho/2-1)$, contradicting deg$(g)=\varrho+1$. Thus $s\ge 3$.
\\

ASSERTION 2: $n=(\varrho+3)/2$.
\\

{\it Proof}. If one had $n>(\varrho+3)/2$, from the above arguments it would follow that 
\[
   {\rm deg}(g)>\frac{\varrho+3}{4}+1+3\biggl(\frac{\varrho+3}{4}-1\biggl).
\]
Hence deg$(g)>\varrho+1$, a contradiction. Therefore $n=(\varrho+3)/2$. \hfill q.e.d.
\ \\

On the one hand, it follows now by (4) that $s=3$ implies $m_3=(1-\varrho)/2\notin\{1,2\}$, a contradiction. Hence $s\ge 4$. On the other hand, $s\ge 6$ gives $\varrho+5\le m_3+m_4+m_5+m_6$. Since $\varrho\ge 4$ and $m_i\le 2$, $\forall$ $i$, this is once again a contradiction. Therefore $4\le s\le 5$. From now on we write $\varrho=4k+1$, $n=2(k+1)$ and so deg$(g)=2(2k+1)$, for $k\in\N^*$. 
\\
  
If $s=5$, then $k=1$ and so $m_3=m_4=m_5=2$, deg$(g)=6$. This means that $M$ has exactly ten points where the normal is vertical. Since ord$(J^2)=2$, then four points of $M$ contribute each with at least $1$ for deg$(g)$, while ord$_{p_i}(g)\ge 3$, $i=1,2$. Hence deg$(g)\ge 4+3>6$, which is absurd. Consequently, $s=4$.
\\

Now (4) simplifies to 
\BE
   \sum_{i=3}^4(2(k+1)-m_i)=4k+2,
\EE
whence $m_3=m_4=1$. Then each fundamental piece of $\tM$ has eight points with vertical normal vectors: two ends $\{p_1,p_2\}$, two points on $Ox_3$ with $m_3=m_4=1$ that we call $\{q_1,q_2\}$, and four points $r_i$, $i=1,\dots,4$, corresponding to $m_1=m_2=2$. Notice that deg$(g)=2(2k+1)$. 
\\

Since $\sigma=J$ and $\tJ$ is a rigid motion, we conclude that $g(r_i)=g(\sigma(r_i))$, $\forall i$. Now it is clear that $g(q_1)=1/g(q_2)$. With no loss of generality we take $r_3=\sigma(r_1)$, $r_4=\sigma(r_2)$, $g(p_1)=g(q_2)=g(r_{1,3})=0$ and $g(p_2)=g(q_1)=g(r_{2,4})=\infty$. Hence, the divisor of $g$ is written as 
\BE
   [g]=\frac{p_1^{k+2}(q_2r_1r_3)^k}{p_2^{k+2}(q_1r_2r_4)^k}.
\EE

Now we are going to write down the divisor of $dh$. For the minimal immersion $X:M\to\R^3/\T$, determined by $(g,dh)$, at each point where $g$ is vertical we must have a zero for $dh$, exactly of the same order as $g$. Moreover, $dh$ must have zeros at the ends $p_{1,2}$ both of order $-2+$ord$(g)_{p_{1,2}}=k$ (see [{\bf 10}, p26] for details). From (6) it follows that  
\BE
   [dh]=(p_1p_2q_1q_2r_1r_2r_3r_4)^k.
\EE

We recall that $\T$ is generated by the vertical translation $\tau(x)=x+(0,0,2)$. So we take a fundamental piece of $\tM$ in the slab $\Ss:=\{(x_1,x_2,x_3)\in\R^3|-1<x_3\le 1\}$ and the points $\tq\in\tM$ such that $\tq/\T\in\{q_1,q_2\}$. Since $J(q_i)=q_i$ and $\tJ$ is a rotation of $\pi/(k+1)$ around $Ox_3$ followed by a reflection in $Ox_1x_2$, then any $\tq$ is in $Ox_3$. Among these points we have $\tq_{1,2}$ in $\Ss$. Therefore, $\tJ(\tq_i)=(0,0,2n_i)-\tq_i$, for some $n_i\in\Z$, $i=1,2$. But $\tJ$ fixes $(0,0,0)$ and $\tJ(\Ss)=\Ss$. We conclude that $\{\tq_1,\tq_2\}=\{(0,0,0),(0,0,1)\}$. Since $J(p_i)=p_i$, we also conclude that the planar ends $\tp_{1,2}$ of $\tM$ are asymptotic to $x_3=0$ and $x_3=1$. We have settled $g(p_1)=0$. Up to changing orientation of $\tM$, $\tp_1$ will correspond to $x_3=1$ and $\tp_2$ to $x_3=0$. 
\\

This means, given a symmetry which fixes one of the points $p_i$, $q_i$, it also fixes the others. Otherwise it interchanges $p_1\leftrightarrow p_2$ and $q_1\leftrightarrow q_2$. We saw already that $\overline{M}/J$ is conformally $\C^2$. Up to a M\"obius transformation one can assume that
\[
   \zeta(p_1)=0,\eh\eh\zeta(p_2)=\infty\eh\eh{\rm and}\eh\eh\zeta(q_1)=1.
\]

Therefore, $\zeta(q_2)$ equals a certain $s\in\C\setminus\{0,1\}$, while $\zeta(r_{1,3})=y_1$ and $\zeta(r_{2,4})=y_2$, namely two distinct complex values in $\C\setminus\{0,1,s\}$. Up to this point, we have not specified the orientation of $\tJ^2$, which can now be fixed as counterclockwise. Let $\gamma_i$ be a single small loop around $0$, $1$, $\infty$, and $y_{1,2}$, for $i=1,\dots,5$, respectively. We take lifts $\hat{\gamma}_i$ of $\gamma_i$ by $\zeta$ and notice that the end points of $\hat{\gamma}_i$ differ by $J^{k_i}$, $0\le k_i\le n-1=2k+1$, $1\le i\le 5$.    
\\

Now take $D$ as the open unitary complex disk at the origin. Since $\zeta$ is the quotient map $(./J):\overline{M}\to\overline{M}/J$, there is a coordinate chart $z:D\to\overline{M}$ with $z(0)=p_1$ such that $\zeta(z)=z^n$. By taking $\gamma_1$ small enough to be in $\zeta(z(D))$, we conclude that $k_1=1$. The same reasoning will give $k_{2,3}=-1$. If we had taken $z(0)=r_1$, then $\zeta(z)=y_1+z^{n/2}$ and so $k_4=2$. By the same reasoning $k_5=-2$. Let us define $\A:=\C\setminus\{0,1,s,y_1,y_2\}$. 
\\

The numbers $k_i$ naturally determine a homomorphism $H:\pi_1(\A)\to\Z_n\oplus\Z_{n/2}$, of which the kernel is $\zeta_*(\pi_1(\overline{M}\setminus\{p_{1,2},q_{1,2},r_{1,..,4}\}))\subset\pi_1(\A)$. Let us now define the following compact Riemann surface:
\BE
   \overline{R}:=\biggl\{(v,w)\in\hat\C\times\hat\C:
   w^{2(k+1)}=\frac{v-s}{v(v-1)}\cdot\biggl(\frac{v-y_1}{v-y_2}\biggl)^2\biggl\}. 
\EE

From (8) one easily sees that $(v,w)\to(v,e^{2\pi i/n}w)$ is a biholomorphism of $\overline{R}$, exactly with the following fixed points: $(0,\infty)$, $(1,\infty)$, $(\infty,0)$, $(s,0)$, each of order $n$, and $(y_1,0)$, $(y_2,\infty)$, each of order $n/2$. The Riemann-Hurwitz formula gives
\[
   \m[4\cdot 1\cdot(n-1)+2\cdot 2\cdot(n/2-1)]-n+1=4k+1, 
\]
namely the same genus as $\overline{M}$. Moreover, the projection map $v:\overline{R}\to\hat\C$, namely $(v,w)\to v$, is such that $v_*(\pi_1(\overline{R}\setminus w^{-1}(\{0,\infty\})))$ also represents the kernel of $H$. From [{\bf 18}, p159] we conclude that $\overline{M}$ is biholomorphic to $\overline{R}$. 
\\

Now we use (6-8) in order to read off the Weierstrass data
\BE
   g=a_0vw^k\eh\eh\eh\eh{\rm and}\eh\eh\eh\eh
   dh=\frac{b_0(v-y_2)w^{k+1}dv}{(v-s)(v-y_1)},
\EE
where $a_0\in\R_{_+}^*$ and $b_0\in((\R_{_+}\times i\R)\setminus(\{0\}\times i\R_{_-}))\subset\C$. These sets containing $a_0$ and $b_0$ were established that way because $\tM$ can be suitably rotated and, if necessary, replaced by its antipodal image. 
\\
\\
\centerline{\bf 4. The symmetries of the minimal immersions}
\

From {\it ({\sl iv})} and the rigidity of $\tM$ we have $|\Delta|\ge 2(\varrho+3)=4n$. Since $|\G|=n$ and $[\Delta:\G]\le 4$, equality follows. Any $A\in\Delta$ either fixes the ends $p_{1,2}$ or interchanges them. In any case $A^2$ fixes the ends, and so $A^2\in\G$. This means that $\Delta/\G=\{f_0,f_1,f_2,f_3\}$ is a group isomorphic to $\Z_2\oplus\Z_2$, and each $f_i$ is an automorphism in $\overline{M}/J\equiv\hat{\C}$. Up to re-indexing we assume that
\begin{itemize}
\item $f_0$ and $f_1$ are holomorphic involutions;
\item $f_2$ and $f_3$ are anti-holomorphic involutions with $f_3=f_1\circ f_2$;
\item $f_0$ and $f_2$ fix the points $0$, $1$, $s$ and $\infty$.
\item $f_1$ and $f_3$ interchange $0\leftrightarrow\infty$ and $1\leftrightarrow s$.
\end{itemize}      

It is immediate that $f_0=id_{\hat\C}$ and $s\in\R\setminus\{0,1\}$, because $f_2$ keeps invariant exactly {\it one} circumference of $\hat\C$, namely $\hat{\R}$. Therefore  
\[
   f_0(v)=v,\eh\eh f_1(v)=\frac{s}{v},\eh\eh f_2(v)=\bar{v}\eh\eh{\rm and}\eh\eh f_3(v)=\frac{s}{\bar{v}}.
\]

We now get more information about $y_{1,2}$. From {\it ({\sl iv})} we have $|Iso(\tM/\T)|\ge 4n$, $n=$ord$(J)$. This implies the existence of automorphisms of $\overline{M}$ that interchange the points $p_1\leftrightarrow p_2$ and $q_1\leftrightarrow q_2$. Recall that a symmetry which fixes one of the points $p_i$, $q_i$ must also fix the others. Hence there exists three distinct automorphisms of $\overline{M}$, $\sigma_{0,1,2}$, such that any $A\in\Delta$ belongs to one of the sets in the following table: 
\[
\begin{tabular}{|c|c|c|c|}\hline
 automorphisms & fix $p_i$    & interch. $p_i$ \\ \hline
 holom.        & $\G$         & $\sigma_0\G$   \\ \hline
 anti-holom.   & $\sigma_1\G$ & $\sigma_2\G$   \\ \hline
\end{tabular}
\]

Let us call $\m\T_\theta$ the translation in $\R^3$ by $(0,0,1)$, followed by a counterclockwise rotation of angle $\theta$ around $Ox_3$. A little reflection about all isometries of $\R^3/\T$ which either fix or interchange $q_{1,2}$ will establish that: 
\begin{itemize}
\item $\sigma_0$ corresponds to a $180^\circ$-rotation around a line $\ell_0\perp Ox_3$ at $x_3=1/2$;
\item $\sigma_1$ can be taken as a $180^\circ$-rotation around a line $\ell_1\perp Ox_3$ at $x_3=0$;
\item $\sigma_2$ corresponds to a screw motion $\m\T_\theta$ for a certain $\theta$.
\end{itemize}   

Notice that $\m\T_\theta$ is orientation-reversing, thus anti-holomorphic. Moreover, $\tM$ can be re-positioned in $\R^3$ in such a way that $\sigma_1(a,b,c)=(a,-b,-c)$. Hence $J\circ\sigma_1(a,b,c)=((a-ib)e^{2\pi i/n},c)$, namely a reflection in the plane $x_2/x_1=\tan(\pi /n)$. This means that we have included reflection in a vertical plane containing $Ox_3$. We have not considered $\sigma_2$ as a reflection in the plane $x_3=1/2$, for this will happen if and only if $\ell_0$ belongs to a vertical plane of reflectional symmetry. 
\\

Now consider the points $\tr_i$ of $\Ss$ such that $\tr_i/\T=r_i$. Therefore, $\sigma_0$ interchanges $r_1\leftrightarrow r_2$ and $r_3\leftrightarrow r_4$, while $\sigma_1$ interchanges $r_1\leftrightarrow r_3$ and $r_2\leftrightarrow r_4$. From $\sigma_1$ we have $y_{1,2}\in\R$ and from $\sigma_0$ it follows that $y_1\cdot y_2=s$.  
\\

Notice that (8) implies 
\[
   v(v-1)(v-s)=\frac{(v-s)^2(v-y_1)^2}{w^{2(k+1)}(v-y_2)^2},
\]
and so $v(v-1)(v-s)$ has a well-defined square root on $\overline{R}$. One rewrites (9) as
\BE
   g=a_0vw^k\eh\eh\eh\eh{\rm and}\eh\eh\eh\eh
   dh=\frac{b_0dv}{\sqrt{v(v-1)(v-s)}}.
\EE

Now we are going to read off some information about the constant $b_0$ at (10). Recall that $\tp_1$ corresponds to $x_3=1$ and $\tp_2$ to $x_3=0$. Therefore, any path on $\tM$ starting at $\tq_i$ and diverging to the end $\tp_1$ will be taken by $\zeta$ to a curve in $\hat{\C}$ connecting $0$ and $\zeta(q_i)$. This latter is homotopically the (oriented) segment $[0,\zeta(q_i))$, for any extra loop with base point at $0$ or $\zeta(q_i)$ gives $Re\oint dh=0$. Our analysis can be now separated in three different cases:

\underline{Case I: $s<0$}. Suppose $(0,0,1)=\tq_1$. By taking $v(t)=t$, $0<t<1$, then $Re\int dh=0\see b_0\in\R_{_+}^*$. On the other hand, if $(0,0,1)=\tq_2$ and $s<t<0$, then $Re\int dh=0\see b_0\in i\R_{_+}^*$.

\underline{Case II: $0<s<1$}. Suppose $(0,0,1)=\tq_1$. In this case, for $\eta:=dt/\sqrt{|t(t-1)(t-s)|}$ one must have 
\BE
   1=Re(b_0)\int_0^s\eta=|Im(b_0)|\int_s^1\eta.
\EE
Therefore, $b_0$ can be neither real nor pure imaginary. Since $[0,1]$ is homotopically $\hat{\R}\setminus(0,1)$, then $1=Re(b_0)\int_1^\infty\eta=|Im(b_0)|\int_{-\infty}^0\eta$, which is equivalent to (11). Indeed, the change $t\to s/t$ shows that $\int_0^s\eta=\int_1^\infty\eta$, while the changes $t\to 1-1/t$ for $\int_{-\infty}^0\eta$ and $t\to 1+(s-1)t$ for $\int_s^1\eta$ show equality between these last two integrals. On the other hand, if $(0,0,1)=\tq_2$ and $0<t<s$, then $Re\int dh=0\see b_0\in i\R_{_+}^*$.

\underline{Case III: $s>1$}. Suppose $(0,0,1)=\tq_1$. By taking $v(t)=t$, $0<t<1$, then $Re\int dh=0\see b_0\in i\R_{_+}^*$. On the other hand, if $(0,0,1)=\tq_2$ and $0<t<s$, one must have 
\BE
   1=Re(b_0)\int_0^1\eta=|Im(b_0)|\int_1^s\eta.
\EE
Therefore, $b_0$ can be neither real nor pure imaginary. Since $[0,s]$ is homotopically $\hat{\R}\setminus(0,s)$, then $1=Re(b_0)\int_s^\infty\eta=|Im(b_0)|\int_{-\infty}^0\eta$, which is equivalent to (12) by suitable changes of variable.
\\

We see that Cases I-III are independent of the real numbers $y_{1,2}$. For any interval $(a,b)\subset\R^*$, since $y_1\cdot y_2=s$, then $y_1\in(a,b)\see y_2\in(s/a,s/b)$ for $s<0$, $y_2\in(s/b,s/a)$ for $s>0$. These are all the possibilities:

\hfil i) $y_{1,2}<s<0<y_{2,1}<1$, \hfil

\hfil ii) $s<y_{1,2}<0<1<y_{2,1}$, \hfil

\hfil iii) $y_{1,2}<y_{2,1}<0<s<1$, \hfil

\hfil iv) $0<y_{1,2}<s<1<y_{2,1}$, \hfil

\hfil v) $0<s<y_{1,2}<y_{2,1}<1$, \hfil

\hfil vi) $y_{1,2}<y_{2,1}<0<1<s$, \hfil

\hfil vii) $0<y_{1,2}<1<s<y_{2,1}$, \hfil

\hfil viii) $0<1<y_{1,2}<y_{2,1}<s$. \hfil

At this point we have just listed the considerable amount of 32 possibilities. Nevertheless, this number will quickly drop to only four items until next section. From (8) and (10) we have
\BE
   g^{2(k+1)}=a_0^{2(k+1)}v^{k+2}\biggl(\frac{v-s}{v-1}\biggl)^k\biggl(\frac{v-y_1}{v-y_2}\biggl)^{2k}. 
\EE

For Case I.i, $(0,0,1)=\tq_1$ implies $b_0>0$. After a suitable rotation of $\tM$ around $Ox_3$, either $v\in(y_1,s)$ or $v\in(y_2,1)$ will give $g\in\R_+^*$. Hence $\phi_2$ is real and never zero on these stretches, and so $Re\int\phi_2\ne 0$. But $\zeta^{-1}(\{1,s,y_{1,2}\})\subset Ox_3$, a contradiction. Therefore I.i implies $(0,0,1)=\tq_2$. For Case I.ii, $(0,0,1)=\tq_2$ implies $ib_0<0$. After a suitable rotation, either $v\in(s,y_1)$ or $v\in(1,y_2)$ will give $g\in\R_+^*$, and the same reasoning leads to the contradiction $Re\int\phi_2\ne 0$. Hence I.ii implies $(0,0,1)=\tq_1$.    
\\

For any of the Cases II.iii, II.iv or II.v, if $(0,0,1)=\tq_1$ then $Re\int dh\ne 0$ for some stretch $v\in(a,b)$, $a,b\in\{s,1,y_{1,2}\}$. Once again, this contradicts $\zeta^{-1}(\{1,s,y_{1,2}\})\subset Ox_3$ and therefore $(0,0,1)=\tq_2$. The same reasoning shows that III.vi, III.vii and III.viii imply $(0,0,1)=\tq_1$.
\\

At this point, since either $y_1<y_2$ or $y_2<y_1$, we have just reduced our analysis to 16 Cases. Before going ahead, notice that both II.iv and III.vii fail. This is because $\tM$ can be suitably rotated about $Ox_3$ to get $g$ real for $\max\{y_1,y_2\}\le v<\infty$, while $dh$ is pure imaginary on this stretch. Therefore, $dh\cdot dg/g$ is pure imaginary there, implying that $\tM/\T$ has horizontal straight lines connecting its ends to points $r_i$, $1\le i\le 4$. Since these points are in $Ox_3$, some of them should coincide with either $q_1$ or $q_2$, contradicting the embeddedness of $\tM$. Now they remain the other 12 cases.       
\\
\\
\centerline{\bf 5. Reduction of cases by geometric arguments}
\

We have just concluded that $(0,0,1)=\tq_2$ exactly for I.i, II.iii and II.v, all with $b_0\in i\R_{_+}^*$, while $(0,0,1)=\tq_1$ exactly for I.ii, III.vi and III.viii, the latter two also with $b_0\in i\R_{_+}^*$, while $b_0\in \R_{_+}^*$ for the former. 
\\

Now consider i, ii, v and viii. If $y_1$ comes immediately after $s$ or vice-versa, then (13) shows that real values of $v$ between $s$ and $y_1$ will make $g$ vary along some meridians of $\hat\C$, from $0$ to $0$, but never reaching $\infty$. Therefore, $dg/g$ will be real, while $dh$ is real and never zero. Hence these stretches are plane geodesics of $\tM$. However, since they connect $\tq_2$ with $\tr_1$ or $\tr_3$, which lie in $Ox_3$, then any of these curves will cross the vertical axis at a third point in between, where the normal vector will not be vertical. But this contradicts the embeddedness of $\tM$, since $\tJ^2$ is a $2\pi/(k+1)$-rotation around $Ox_3$. Therefore, cases i, ii, v and viii are reduced to
\\

\hfil i) $y_2<s<0<y_1<1$, \hfil

\hfil ii) $s<y_2<0<1<y_1$, \hfil

\hfil v) $0<s<y_2<y_1<1$, \hfil

\hfil viii) $0<1<y_1<y_2<s$. \hfil
\\

In Section 3, we concluded that $\tJ$ is a rotation around $Ox_3$ followed by reflection in $Ox_1x_2$. Therefore, up to re-indexing, we assume that $r_{1,2}$ lie between $x_3=0$ and $x_3=1$. If $r_1$ is above $r_2$, this will force $y_2<y_1$ at iii and vi due to $Re\int dh$, namely the third coordinate of $\tM$. But $y_2<y_1$ drops case iii, for (13) shows that this would give two geodesics in the plane $\Pi:x_2/x_1=-\tan\frac{k\pi}{2(k+1)}$, the one bounded and connecting $q_1$ with $q_2$, the other unbounded and connecting the end $p_2$ with $r_2$. Moreover, by (13) ones sees that both should be in the same half-plane determined by $Ox_3\subset\Pi$ and therefore would cross, contradicting the embeddedness of $\tM$.    
\\

The assumption of $r_1$ above $r_2$ also drops case ii because of the following argument: after a suitable rotation of $\tM$ around $Ox_3$, $v\in(y_2,0)$ will give $g\in\R_+^*$. Hence, the corresponding geodesic in $\tM$ will have to cross the vertical axis at a point where the normal vector will not be vertical. But this contradicts the embeddedness of $\tM$, since $\tJ^2$ is a $2\pi/(k+1)$-rotation about $Ox_3$. This means, $r_1$ above $r_2$ cancels cases ii and iii. By very similar arguments, assuming $r_2$ above $r_1$, cases ii and iii fail again. Now we remain with 

\hfil i) $y_2<s<0<y_1<1$, \hfil

\hfil v) $0<s<y_2<y_1<1$, \hfil

\hfil vi) $y_2<y_1<0<1<s$, \hfil

\hfil viii) $0<1<y_1<y_2<s$. \hfil
\\

Notice from statement vi that $r_1$ is above $r_2$ if and only if $y_2<y_1$. It is totally equivalent to study this one or its reverse, and hence our reduction is complete. The next sections are devoted to the study of the remaining cases.
\\
\\
\centerline{\bf 6. Non-solvability of the period problems for i, v and vi}
\

We begin with case i. In Section 3 one saw that $\overline{M}$ is biholomorphic to $\overline{R}$ given by (8). Now we closely follow the arguments from [{\bf 16}, p453]. By labelling $a=1/\sqrt{|s|}$, $b=|y_1/y_2|^\m$, $B=\sqrt[k+1]{b}$ and making the changes $z=av$, $u=B/w$, an easy computation shows that $\overline{R}$ is biholomorphic to
\[
   \overline{N}:=\biggl\{(z,u)\in\hat\C\times\hat\C:
   u^{2(k+1)}=\frac{z(z-a)}{az+1}\cdot\biggl(\frac{bz+1}{z-b}\biggl)^2\biggl\}. 
\]   

Notice that we still have $J(z,u)=(z,e^{i\pi/(k+1)}u)$. Define $A:=a_0B^k/a$, so with $u$ and $z$ the Weierstrass data become
\[
   g=Azu^{-k}\eh\eh\eh\eh{\rm and}\eh\eh\eh\eh
   dh=\frac{ab_0dz}{\sqrt{z(z-a)(az+1)}}.
\]

In Section 4 we saw that $\overline{M}$ is endowed with an automorphism $\sigma_2$, which in $\tM$ corresponds to a screw motion $\m\T_\theta$, for a certain $\theta$. From the Weierstrass data, it is easy to see that $\theta=0$ and so $\sigma_2$ represents half of the vertical translation $\T$. Now recall that $b_0\in i\R_+^*$. From arguments very similar to [{\bf 16}, p461], one concludes that
\[
   \sigma_2(z,u)=(-1/\bar{z},1/\bar{u})
\] 
and
\BE
   \sigma_2(gdh)=-A^2\overline{dh/g},\eh\eh\sigma_2(dh/g)=-\overline{gdh}/A^2
   \eh\eh{\rm and}\eh\eh\sigma_2(dh)=\overline{dh}.
\EE

As $\sigma_2$ corresponds to a rigid motion in $\R^3$, then $A^2=1$. Hence $A=1$ because $A$ is positive. Defining $N:=\overline{N}\setminus z^{-1}(\{0,\infty\})$, the immersion $X:N\to\R^3/\T$ will be period free if and only if
\BE
   \int_\gamma dh/g=\int_\gamma\overline{gdh},\eh\eh Re\int_\gamma dh\in 2\Z,
   \eh\eh\forall\eh[\gamma]\in H_1(N).
\EE

\input epsf
\begin{figure}[ht]
\centerline{
\epsfxsize 13cm
\epsfbox{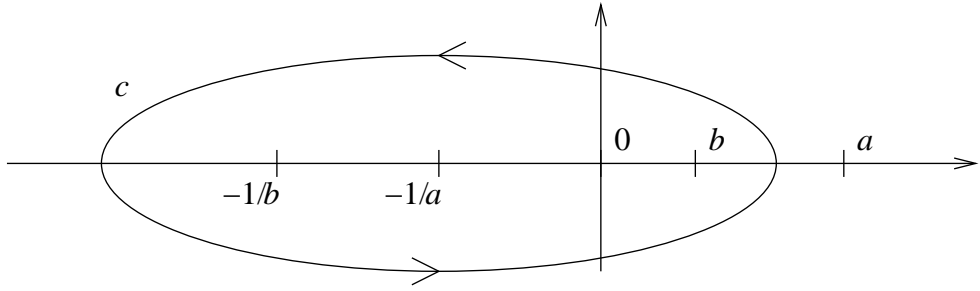}}
\caption{The curve $c$ in the complex plane.}
\end{figure}

Now observe that $b=ay_1$, with $0<y_1<1$, and consider the curve $c$ represented in Figure 3. For $\tg:=e^{-i\pi/(2k+2)}g$, a suitable choice of logarithmic branch shows that the Weierstrass data $(\tg,dh)$ take the intervals $a<z<b$ and $-1/b<z<-1/a$ to geodesics in planes parallel to $x_1=0$. Let $c^+$ denote the stretch of $c$ in the upper half plane. If $X$ is period free, then
\[
   Re\int_{c^+}\tg dh=Re\int_{c^+}dh/\tg,
\]
or equivalently
\BE
   -Re\int_{-1/b<z<b}\tg dh=Re\int_{a<z<\infty}dh/\tg.
\EE

From (14) we see that 
\[
   Re\int_{a<z<\infty}dh/\tg=Re\int_{a<z<\infty}\sigma_2(\tg dh)=-Re\int_{-1/a<z<0}\tg dh.
\]
Now (16) becomes
\BE
   Re\int_{-1/b<z<-1/a}\tg dh+Re\int_{0<z<b}\tg dh=0.
\EE

By keeping the same logarithmic branch, on $-1/b<z<-1/a$ we see that $dh$ is positive and $\tg\in e^{i\pi(k-1)/(2k+2)}\R_+$, while on $0<z<b$ on has $dh$ negative and $\tg\in e^{i\pi(k-1)/(2k+2)}\R_-$. Hence (17) is equivalent to 
\BE
   \cos\biggl(\frac{\pi}{2}\cdot\frac{k-1}{k+1}\biggl)
   \cdot\biggl[\int_{-1/b}^{-1/a}|gdh|+\int_0^b|gdh|\biggl]=0,
\EE
which never holds. 
\\

A three-dimensional sketch of $X(N)$ is presented in Figure 4. It is important to note that the arguments presented herein differ from the ones in [{\bf 16}, pp460-2] and [{\bf 2}, pp176-180]. Of course, the surfaces are not the same, but the arguments are adaptable. For instance, one could rewrite (18) in [{\bf 16}, p462] by suppressing the first integral and taking $b=a$ in the second.  
\\
\input epsf
\begin{figure}[ht]
\centerline{
\epsfxsize 15cm
\epsfbox{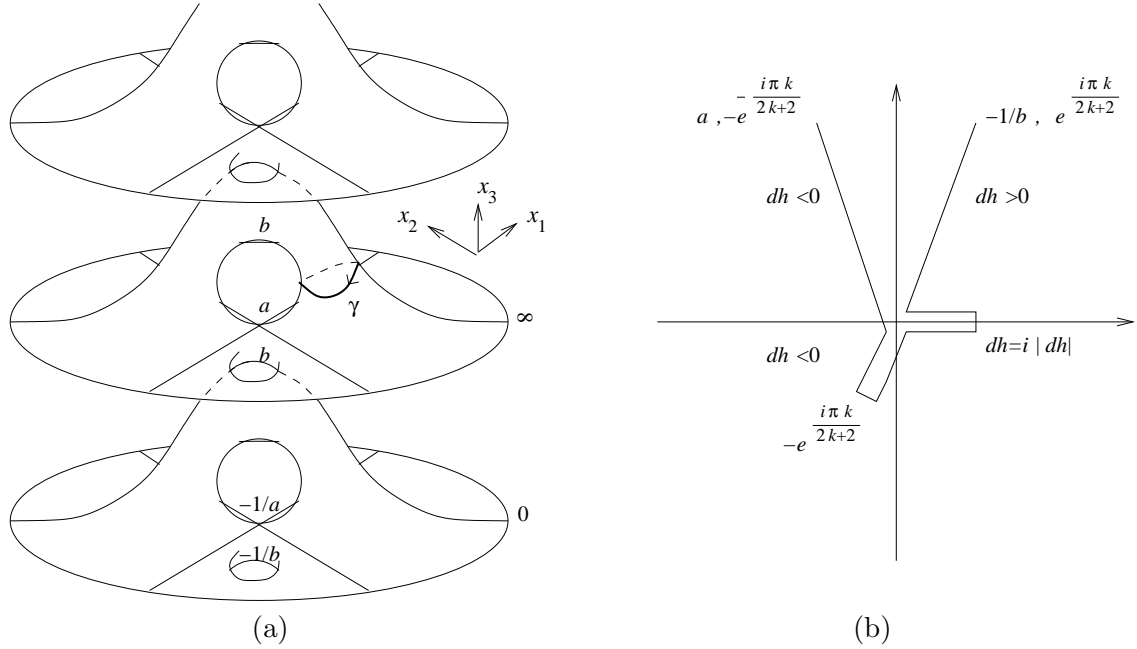}}
\hspace{3.5cm}(a)\hspace{7.5cm}(b)
\caption{(a) The surface $X(N)$ for $k=1$; (b) The Gau\ss \ map on symmetry curves.}
\end{figure}

For reasons that will soon be clear, we shall invert order and study case vi before v. Take $a$, $b$, $B$, $z$ and $u$ as before. Hence $\overline{R}$ is biholomorphic to
\[
   \overline{\Nn}:=\biggl\{(z,u)\in\hat\C\times\hat\C:
   u^{2(k+1)}=\frac{z(z-a)}{az-1}\cdot\biggl(\frac{bz+1}{z+b}\biggl)^2\biggl\}, 
\]   
and we still have $J(z,u)=(z,e^{i\pi/(k+1)}u)$. Take $A$ as before, so with $u$ and $z$ the Weierstrass data become
\BE
   g=Azu^{-k}\eh\eh\eh\eh{\rm and}\eh\eh\eh\eh
   dh=\frac{ab_0dz}{\sqrt{z(z-a)(az-1)}}.
\EE

From Section 4, the automorphism $\sigma_2$ corresponds to a screw motion $\m\T_\theta$, now with $\theta=-\pi/(k+1)$. Again $b_0\in i\R_+^*$, hence 
\[
   \sigma_2(z,u)=(1/\bar{z},e^{\frac{-i\pi}{k+1}}/\bar{u})
\] 
and
\BE
   \sigma_2(gdh)=e^{-\frac{i\pi}{k+1}}A^2\overline{dh/g},\eh\eh
   \sigma_2(dh/g)=e^{\frac{i\pi}{k+1}}\overline{gdh}/A^2
   \eh\eh{\rm and}\eh\eh\sigma_2(dh)=\overline{dh}.
\EE

As $\sigma_2$ corresponds to a rigid motion in $\R^3$, then $A=1$. Defining $\Nn:=\overline{\Nn}\setminus z^{-1}(\{0,\infty\})$, the immersion $X:\Nn\to\R^3/\T$ will be period free if and only if
\[
   \int_\gamma dh/g=\int_\gamma\overline{gdh},\eh\eh Re\int_\gamma dh\in 2\Z,
   \eh\eh\forall\eh[\gamma]\in H_1(\Nn).
\]

\input epsf
\begin{figure}[ht]
\centerline{
\epsfxsize 13cm
\epsfbox{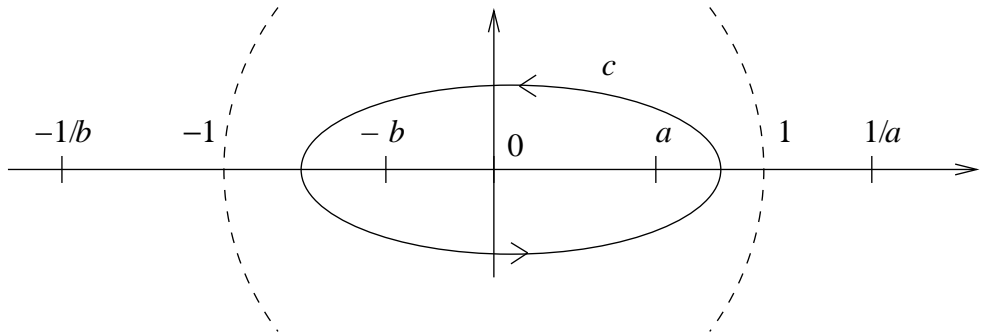}}
\caption{The curve $c$ in the complex plane.}
\end{figure}

Now observe that $a$ and $b$ are both in $(0,1)$, and consider the curve $c$ represented in Figure 5. Take $\tg$ as before. A suitable choice of logarithmic branch shows that the Weierstrass data $(\tg,dh)$ take the intervals $a<z<1/a$ and $0<z<\infty$ to geodesics in planes parallel to $x_1=0$. Let $c^+$ denote the stretch of $c$ in the upper half plane. If $X$ is period free, then
\[
   Re\int_{c^+}\tg dh=Re\int_{c^+}dh/\tg,
\]
or equivalently
\BE
   -Re\int_{-b<z<a}\tg dh=Re\int_{1/a<z<-1/b}dh/\tg,
\EE
with $z$ passing through infinity. The right-hand side of (21) equals 
\[
   Re\int_{1/a<z<-1/b}\overline{dh/\tg}\stackrel{(20)}{=}
   Re\int_{1/a<z<-1/b}e^{\frac{i\pi}{k+1}}\sigma_2(\tg dh)=
   Re\int_{-b<z<a}e^{\frac{i\pi}{k+1}}\tg dh.
\]

Hence (21) holds if and only if $2\cos(\frac{\pi}{2k+2})\cdot Re\int_{-b}^agdh=0$. But since $gdh$ is pure imaginary for $0<z<a$, (21) is equivalent to $Re\int_{-b}^0gdh=0$. This never holds, for
\[
   Re\int_{-b}^0gdh=-Re\int_{-b}^0e^{\frac{i\pi}{2k+2}}|gdh|=
   -\cos\biggl(\frac{\pi}{2k+2}\biggl)\cdot\int_{-b}^0|gdh|<0.
\]

A three-dimensional sketch of $X(\Nn)$ is depicted in Figure 6. As remarked at the introduction, one easily identifies the presence a Gaussian geodesic in cases i and vi. One period is essentially due to $\int|gdh|$ along it, and consequently never vanishes. That is why we surveyed cases i and vi together. In case v there is no Gaussian geodesic, but Horgan saddle. The reader will notice that one period never vanishes again, namely around that saddle. In spite of their oddness, Horgan saddles were recently found in singly periodic examples (see [{\bf 30}]). 
\\
\input epsf
\begin{figure}[ht]
\centerline{
\epsfxsize 15cm
\epsfbox{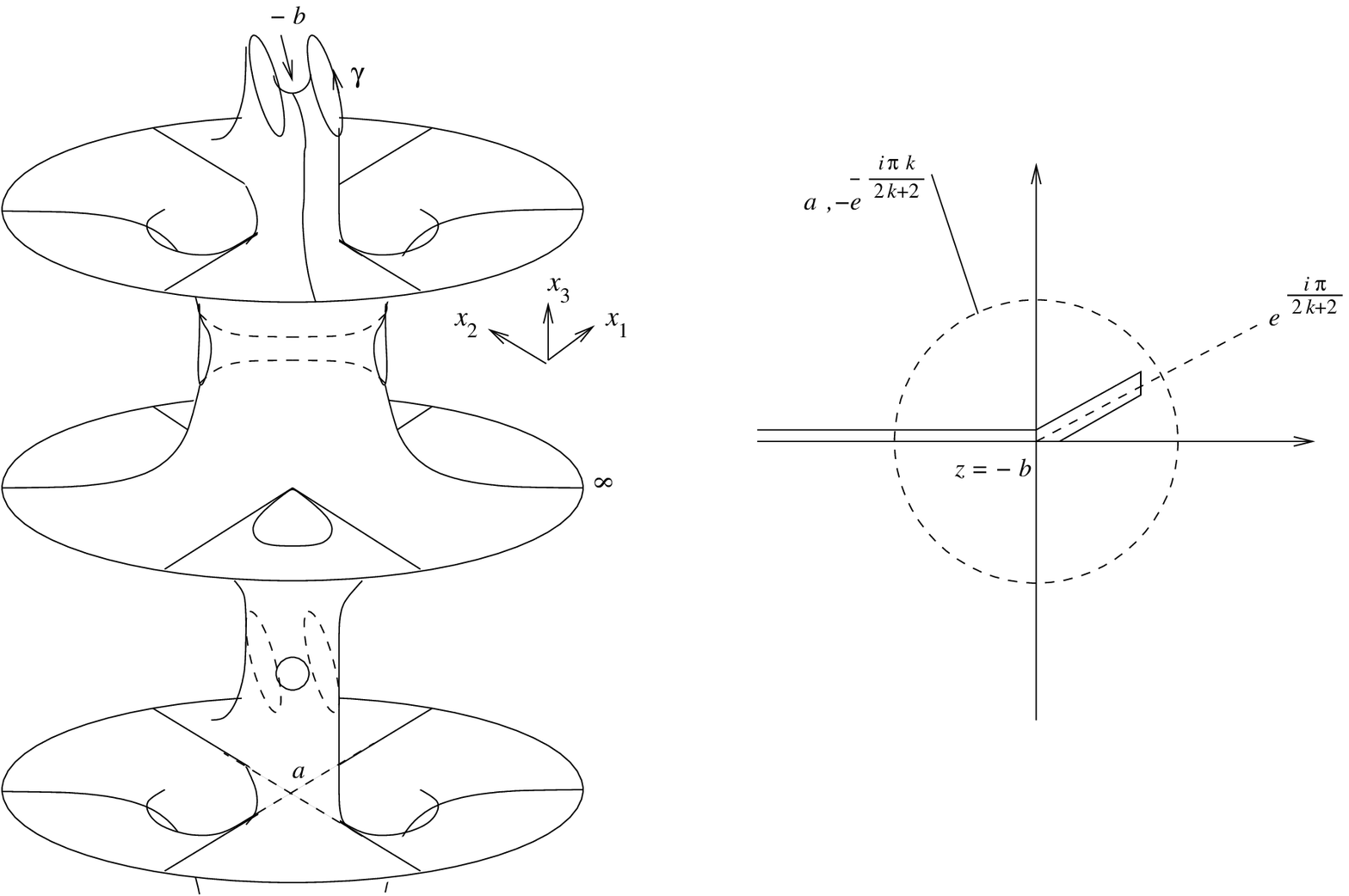}}
\hspace{3.5cm}(a)\hspace{7.5cm}(b)
\caption{(a) The surface $X(\Nn)$ for $k=1$; (b) The Gau\ss \ map on symmetry curves.}
\end{figure}

We conclude this section with case v. By taking all parameters as before, one sees that $\overline{R}$ is biholomorphic to  
\[
   \overline{\Ns}:=\biggl\{(z,u)\in\hat\C\times\hat\C:
   u^{2(k+1)}=\frac{z(z-a)}{az-1}\cdot\biggl(\frac{bz-1}{z-b}\biggl)^2\biggl\}, 
\]   
again with $J(z,u)=(z,e^{i\pi/(k+1)}u)$. The Weierstrass data coincide with (19), but now $a\in(1,\infty)$ and $b\in(1,a)$. The automorphism $\sigma_2$ is once more the screw motion $\m\T_\theta$ with $\theta=\pi/(k+1)$ and (20) still holds, with $A=1$. Defining $\Ns:=\overline{\Ns}\setminus z^{-1}(\{0,\infty\})$, the immersion $X:\Ns\to\R^3/\T$ will be period free if and only if
\[
   \int_\gamma dh/g=\int_\gamma\overline{gdh},\eh\eh Re\int_\gamma dh\in 2\Z,
   \eh\eh\forall\eh[\gamma]\in H_1(\Ns).
\]

Now consider the curve $c$ represented in Figure 7. By using $\tg$ and a suitable choice of logarithmic branch, the pair $(\tg,dh)$ takes the intervals $b<z<a$ and $1/a<z<1/b$ to geodesics in planes parallel to $x_1=0$.
\\

\input epsf
\begin{figure}[ht]
\centerline{
\epsfxsize 13cm
\epsfbox{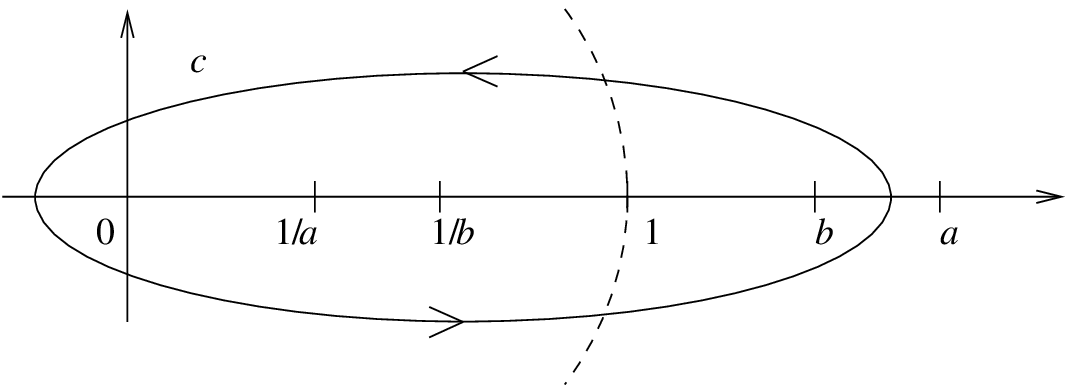}}
\caption{The curve $c$ in the complex plane.}
\end{figure}

Let $c^+$ denote the stretch of $c$ in the upper half plane. If $X$ is period free, then
\[
   Re\int_{c^+}\tg dh=Re\int_{c^+}dh/\tg,
\]
or equivalently
\BE
   -Re\int_{0<z<b}\tg dh=Re\int_{a<z<\infty}dh/\tg.
\EE

From (20) we see that 
\[
   Re\int_{a<z<\infty}dh/\tg=Re\int_{a<z<\infty}e^{\frac{i\pi}{k+1}}\sigma_2(\tg dh)=
   Re\int_{0<z<1/a}e^{\frac{i\pi}{k+1}}\tg dh,
\]
whereas 
\BE
   Re\int_{0<z<b}\tg dh=Re\int_{0<z<1/a}\tg dh+Re\int_{1/b<z<b}\tg dh,
\EE
because $\tg$ and $dh$ are real and pure imaginary for $1/a<z<1/b$, respectively. For $0<z<1/a$, $g=-|g|$ and $dh=i|dh|$, and so the last integral from (23) cancels with the right-hand side of (22). This leads to $Re\int_{0<z<1/a}\tg dh=0$, which never holds. See Figure 8 for a sketch of $X(\Ns)$.
\input epsf
\begin{figure}[ht]
\centerline{
\epsfxsize 15cm
\epsfbox{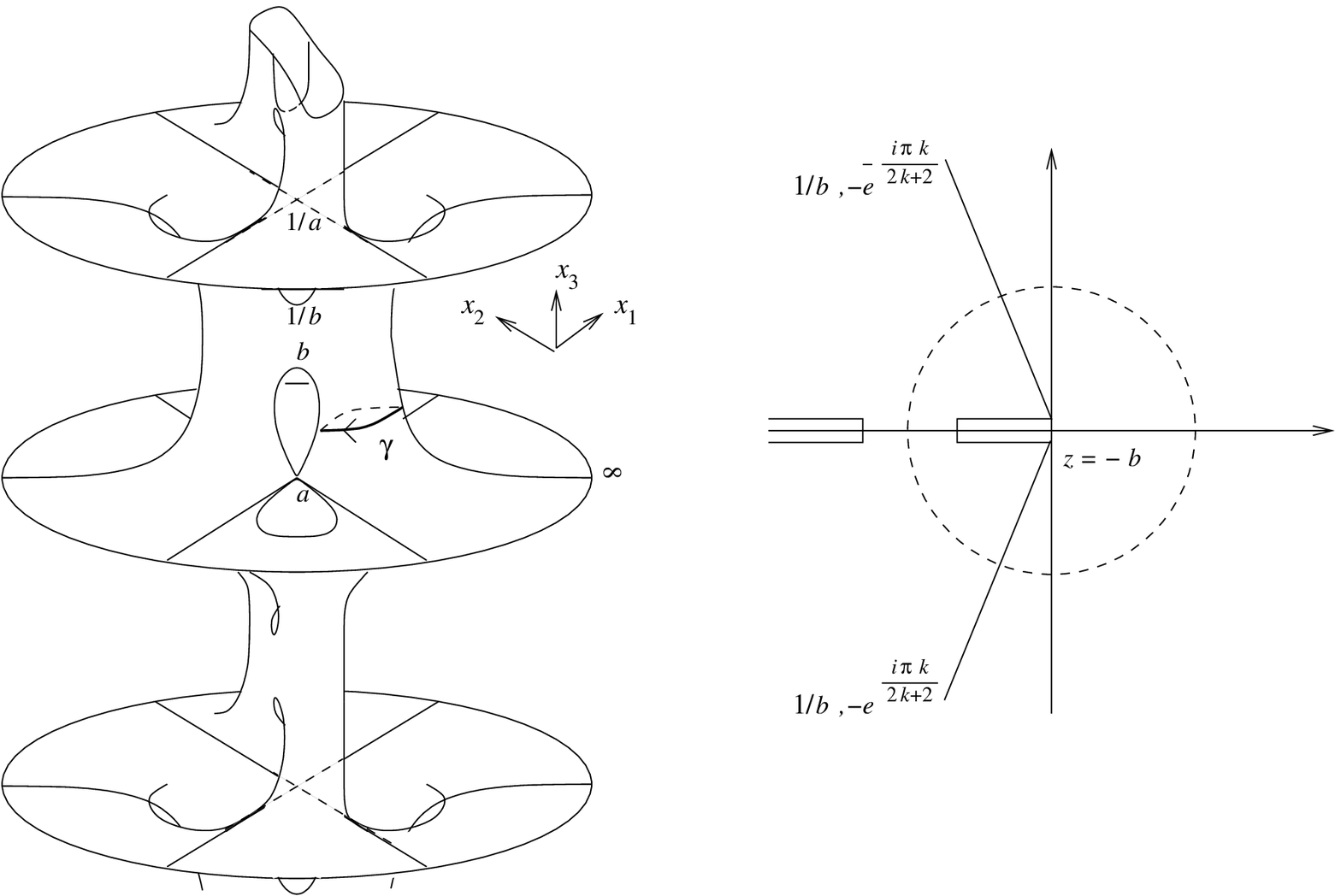}}
\hspace{3.5cm}(a)\hspace{7.5cm}(b)
\caption{(a) The surface $X(\Ns)$ for $k=1$; (b) The Gau\ss \ map on symmetry curves.}
\end{figure}
\ \\
\\
\centerline{\bf 7. The Hoffman-Wohlgemuth surfaces}
\

We recall that $0<1<y_1<y_2<s$ for case viii, thus with all parameters positive as in case v. Therefore $\overline{R}$ is again biholomorphic to $\overline{\Ns}$, but now with $a\in(0,1)$ and $b\in(a,1)$. Both (19) and (20) still hold, hence $A=1$. 
\\

Now consider the curves $\gamma$, $\Gamma$, $\delta$ and $\Delta$ represented in Figure 9. Up to homotopy, the curve $\delta+\Delta=\Gamma-\gamma$ is invariant under the map $z\to 1/z$. 
\\

\input epsf
\begin{figure}[ht]
\centerline{
\epsfxsize 13cm
\epsfbox{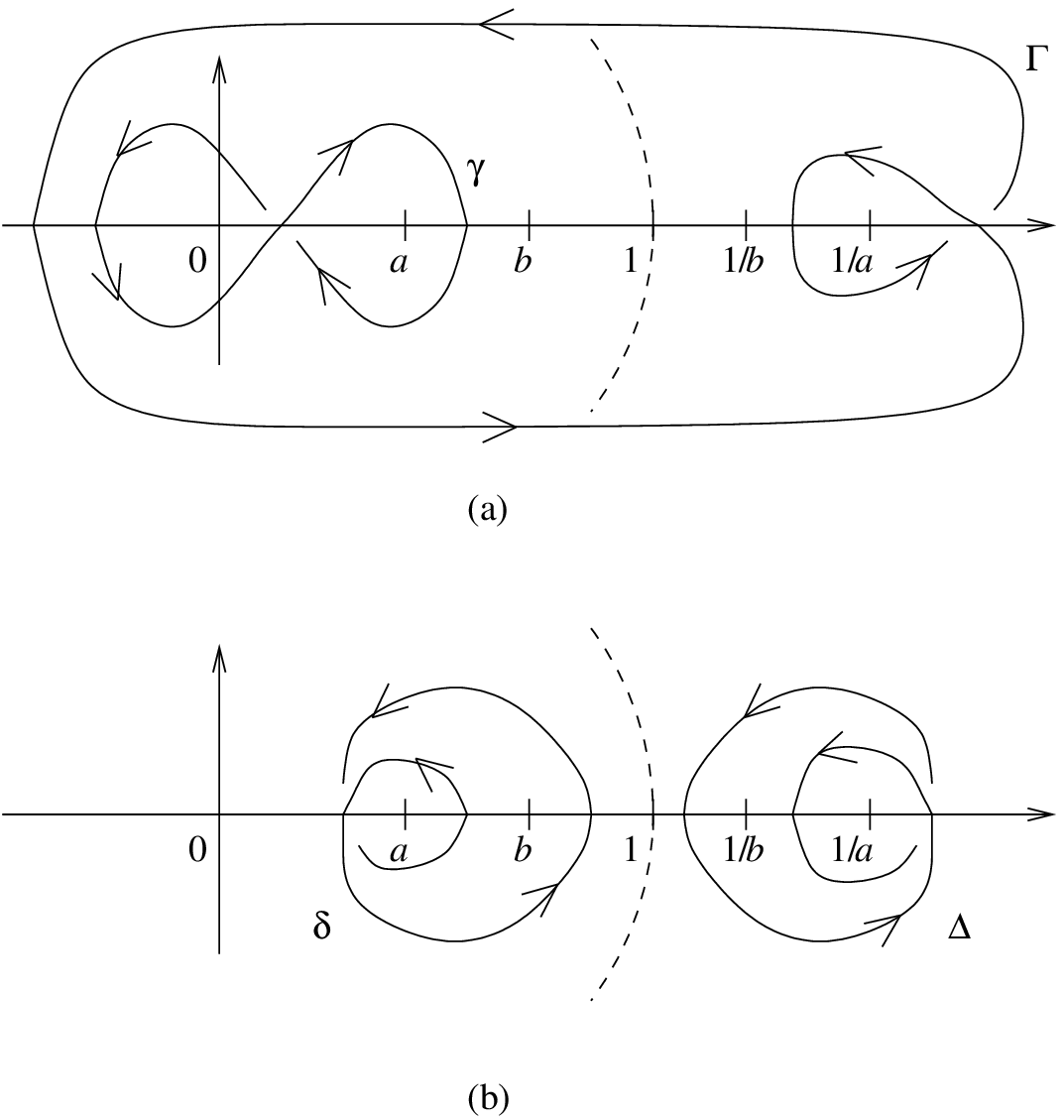}}
\caption{(a) The curves $\gamma$, $\Gamma$ and (b) $\delta$, $\Delta$ in the complex plane.}
\end{figure}

Recalling that $\tg=e^{-i\pi/(2k+2)}g$ and defining $\hat{g}:=e^{-i\pi k/(2k+2)}g$, the immersion $X:\Ns\to\R^3/\T$ will be period free if and only if
\BE
   Re\int_{\gamma}\tg dh=Re\int_{\gamma}dh/\tg
\EE
and
\BE
   Re\int_{\delta}\hat{g}dh=Re\int_{\delta}dh/\hat{g}.
\EE 

Let us take $\psi:=\sigma_0\circ\sigma_1\circ\sigma$, hence $\psi(z,u)=(1/\bar{z},1/\bar{u})$, $\psi(g)=1/\bar{g}$ and $\psi(dh)=-\overline{dh}$. Consequently, 
\BE
   \int_{\delta+\Delta}gdh=-\int_{\psi(\delta+\Delta)}gdh=-\int_{\delta+\Delta}\psi(gdh)=
   \int_{\delta+\Delta}\overline{dh/g}.
\EE

Moreover, $Re\int_{\gamma}dh/\tg=Re\int_{\Gamma}\overline{dh/\tg}-Re\int_{\delta+\Delta}\overline{dh/\tg}\stackrel{(26)}{=}-Re\int_{\Gamma}\psi(\tg dh)-Re\int_{\delta+\Delta}\tg dh$, namely
\BE
   Re\int_{\gamma}dh/\tg=-Re\int_{\gamma}\tg dh-Re\int_{\delta+\Delta}\tg dh.
\EE

With (27) one rewrites (24) as
\BE
   2Re\int_{\gamma}\tg dh=-Re\int_{\delta+\Delta}\tg dh.
\EE

Up to a homothety in $\R^3$, appropriate choices of the logarithmic branch will give
\[
   \tg(t)=e^{-\frac{i\pi}{2k+2}}t^\frac{k+2}{2k+2}\cdot
          \biggl(\frac{1-at}{a-t}\biggl)^\frac{k}{2k+2}
          \cdot\biggl(\frac{b-t}{1-bt}\biggl)^\frac{k}{k+1},
   dh(t)=\frac{idt}{\sqrt{t(a-t)(1-at)}},\eh 0<t<a,
\]
\[
   \tg(t)=e^{\frac{i\pi(k-1)}{2k+2}}t^\frac{k+2}{2k+2}\cdot
          \biggl(\frac{1-at}{t-a}\biggl)^\frac{k}{2k+2}
          \cdot\biggl(\frac{b-t}{1-bt}\biggl)^\frac{k}{k+1},
   dh(t)=\frac{dt}{\sqrt{t(t-a)(1-at)}},\eh t^{\pm 1}\in(a,b),
\]
and
\[
   \hat{g}(t)=t^\frac{k+2}{2k+2}\cdot\biggl(\frac{1-at}{t-a}\biggl)^\frac{k}{2k+2}
          \cdot\biggl(\frac{b-t}{1-bt}\biggl)^\frac{k}{k+1},
   dh(t)=\frac{dt}{\sqrt{t(t-a)(1-at)}},\eh a<t<b.
\] 

It is easy to check that $z\to 1/\bar{z}$ is now an isometry for the minimal surface, which implies $\int gdh=0$ on $b<z<1/b$. Therefore (25) and (28) are explicitly given by
\BE
   I_0:=\int_a^b\frac{t^{1/(2k+2)}[(b-t)/(1-bt)]^{k/(k+1)}dt}{[(1-at)(t-a)^{2k+1}]^{1/(2k+2)}}=
   \int_a^b\frac{[(1-bt)/(b-t)]^{k/(k+1)}dt/t}{[t(t-a)(1-at)^{2k+1}]^{1/(2k+2)}}=:I_1  
\EE
and 
\BE
   J_0:=\int_0^a\frac{t^{1/(2k+2)}[(b-t)/(1-bt)]^{k/(k+1)}dt}{[(1-at)(a-t)^{2k+1}]^{1/(2k+2)}}=
   \cos\biggl(\frac{\pi}{2k+2}\biggl)\cdot J_1, 
\EE
where $J_1:=J_++J_-$ with
\BE 
   J_\pm:=\int_{(a,b)^{\pm 1}}
   \frac{t^{1/(2k+2)}[(b-t)/(1-bt)]^{k/(k+1)}dt}{[(1-at)(t-a)^{2k+1}]^{1/(2k+2)}}.
\EE

Notice that $I_0=J_+$ and $I_1=J_-$. Except for $I_1$ and $J_-$, it is easy to see that all integrals in (29-31) are continuous at $b=1$. For $I_1$ and $J_-$ make the changes $t=b-u^{k+1}$ and $t=(1+u^{k+1})/b$, respectively. Take the limit $b\to 1$ and make back the change $u=\sqrt[k+1]{|1-t|}$. The functions ``$f_1(a)$'' and ``$f_2(a)$'' described in [{\bf 16}, p457] are exactly $J_1(a,b)$ and $J_0(a,b)$ at $b=1$, respectively. From this point on we shall follow some ideas from [{\bf 16}, p457-9].
\\

The change $t=(1/a-a)u+a$ gives
\[
   J_1(a,1)=a^{\frac{1}{k+1}}
   \int_0^1\biggl(\frac{a^{-2}+(a^{-4}-a^{-2})t}{(1-t)t^{2k+1}}\biggl)^\frac{1}{2k+2}dt,
\]
and so 
\BE
   \frac{\deh J_1(a,1)}{\deh a}=\frac{J_1(a,1)}{a(k+1)}+{\rm (neg.term)}.
\EE

Moreover, the change $t=au$ for $J_0(a,1)$ will give
\BE
   \frac{\deh J_0(a,1)}{\deh a}=\frac{J_0(a,1)}{a(k+1)}+{\rm (pos.term)}.
\EE

Combining (32) with (33) it follows that $J_0(a,1)/J_1(a,1)$ is strictly increasing. By taking $\B$ as the beta function, we now closely follow the computations from [{\bf 16}, pp457-8] to conclude that
\[
   \lim_{a\to 0}(a^\frac{1}{k+1}J_+(a,b))=0,
\]
\[
   \lim_{a\to 0}(a^\frac{1}{k+1}J_-(a,b))=
   b^\frac{-k}{k+1}\B\biggl(\frac{1}{k+1},\frac{2k+1}{2k+2}\biggl),
\]
and
\[
   \lim_{a\to 0}(a^\frac{-1}{k+1}J_0(a,b))=
   b^\frac{k}{k+1}\B\biggl(\frac{1}{2k+2},\frac{2k+3}{2k+2}\biggl).
\]

Namely, for $a$ close to zero one has $J_0<\cos(\pi/(2k+2))\cdot J_1$. Up to this point, we have been considering $(a,b)\in(0,1)\times(a,1)$. However, the equivalent choice $(b,a)\in(0,1)\times(0,b)$ will be easier to deal with. Now notice that
\[
   J_+\le({\rm pos.const.})\int_a^b\frac{(b-t)^{k/(k+1)}dt}{(t-a)^{(2k+1)/(2k+2)}},
\]
and so the change $t=(b-a)u+a$ gives
\[
   J_+\le({\rm pos.const.})\int_0^1
   \frac{(b-a)^{(2k+1)/(2k+2)}(1-u)^{k/(k+1)}du}{u^{(2k+1)/(2k+2)}}\to 0,\eh{\rm for}\eh a\to b.
\]
Moreover,
\[
   J_-\le({\rm pos.const.})\int_{1/b}^{1/a}\frac{(1-at)^{-1/(2k+2)}dt}{(bt-1)^{k/(k+1)}},
\]
and so the change $t=(1/a-1/b)u+1/b$ gives 
\[
   J_-\le({\rm pos.const.})\int_0^1
   \frac{(1/a-1/b)^{1/(2k+2)}du}{(1-u)^{1/(2k+2)}u^{k/(k+1)}}\to 0,\eh{\rm for}\eh a\to b.
\]

On the one hand, if $a\to b$ then $J_\pm$ will both vanish. On the other hand, $J_0$ will remain finite and positive. Therefore $J_0>\cos(\pi/(2k+2))\cdot J_1$ when $a$ is close to $b$. Let us define $\D:=\{z\in\C:\eh 0<Im\eh z<\eh Re\eh z<1\}$. Until this point we have that
\\

\it There is an analytic curve $\af:(0,1)\to\D$ for which any $s\in(0,1)$ will make the choice $(b,a)=\af(s)$ a solution for (30). Up to orientation reversing, $\lim_{s\to 0}\af(s)=(0,0)$ and $\lim_{s\to 1}\af(s)=(1,a_1)$, for some $a_1\in(0,1)$. Moreover, when $s$ is sufficiently close to 1, $\af(s)$ is a graph of\eh\eh $a$ as function of\eh\eh $b$.\rm   
\\

From now on we shall always work with $(b,a)=\af(s)$, for some $s\in(0,1)$. For the integrands of (29), make the change $t=\eps u+a$ with $\eps=b-a$. One easily computes
\BE
   \lim_{s\to 0}\frac{I_0}{\eps^\frac{1}{2k+2}}=
   \lim_{s\to 0}\eps^\frac{2k}{2k+2}\int_0^1
   \biggl(\frac{t(u)}{(1-at(u))u^{2k+1}}\biggl)^{1/(2k+2)}\cdot
   \biggl(\frac{1-u}{1-bt(u)}\biggl)^{k/(k+1)}du=0
\EE   
and
\BE
   \lim_{s\to 0}\frac{I_1}{\eps^\frac{1}{2k+2}}=
   \lim_{s\to 0}\int_0^1
   \frac{[(1-bt(u))/(1-u)]^{k/(k+1)}}{[u(1-at(u))^{2k+1}t(u)^{2k+3}]^{1/(2k+2)}}du=\infty.
\EE 

The careful reader must have noticed that, intuitively, the extreme case $b=1$ corresponds to the Callahan-Hoffman-Meeks surfaces $M_k$, described in [{\bf 3}]. Their underlying Riemann surfaces have lower genera, and so it is hard to formalise any convergence statement. However, the Callahan-Hoffman-Meeks surfaces were again described in [{\bf 16}], where the integrals for the period problem coincide with $J_0|_{b=1}$ and $J_1|_{b=1}$. In the case of [{\bf 16}], the integration of $Re(1/\hat{g}-\hat{g},i/\hat{g}+i\hat{g},2)dh$ along $\delta+\Delta$ is a geodesic in the plane $x_2=0$, connecting the saddles $z=a$ and $z=1/a$, and symmetric under reflection in the plane $x_3=1/2$. In [{\bf 3}] those surfaces were proved to be embedded, and therefore the geodesic cannot cross the vertical axis, except at $z=a$ and $z=1/a$. Hence $Re\int_\delta(1/\hat{g}-\hat{g})dh=I_1|_{b=1}-I_0|_{b=1}<0$. Together with (34) and (35), this gives $s^*\in(0,1)$ for which $(b,a)=\af(s^*)$ simultaneously solves (29) and (30).      
\\
\\
\centerline{\bf 8. Embeddedness of the Hoffman-Wohlgemuth surfaces}
\

In this last section we use arguments very similar to [{\bf 10}, p60-2] or [{\bf 26}, p360-2]. For convenience of the reader, we recall that every $k\in\N^*$ admits a well-defined complete minimal immersion $X:\Ns\to\R^3/\T$, where $\Ns=\overline{\Ns}\setminus z^{-1}(\{0,\infty\})$,   
\[
   \overline{\Ns}=\biggl\{(z,u)\in\hat\C\times\hat\C:
   u^{2(k+1)}=\frac{z(z-a)}{az-1}\cdot\biggl(\frac{bz-1}{z-b}\biggl)^2\biggl\}, 
\]   
and $X$ is given by the Weierstrass pair
\[
   g=\frac{z}{u^k}\eh\eh\eh\eh{\rm and}\eh\eh\eh\eh
   dh=\frac{ab_0dz}{\sqrt{z(z-a)(az-1)}}.
\]

Here, $a$ and $b$ were determined in the previous section, and $b_0\in i\R_+^*$ is such that $\int_{-\infty}^0 dh=1$. Now consider the domain $\Dd:=\{z\in\C:|z|<1<1+Im(z)\}$. By choosing appropriate logarithm branches, one sees that $g(\Dd)$ is contained in a hemisphere of $\hat\C$ (see Figure 10). As before, take $\tg=e^{-i\pi/(2k+2)}g$ and call $\tX:\Ns\to\R^3/\T$ the minimal immersion given by $(\tg,dh)$. If $(x_1,x_2,x_3)$ are the coordinates of $\tX$, then $(x_2,x_3):\Dd\to\R^2$ is an immersion, and so its image boundary coincides with $(x_2,x_3)(\deh\Dd\setminus\{0\})$.
\\    
\input epsf
\begin{figure}[ht]
\centerline{
\epsfxsize 15cm
\epsfbox{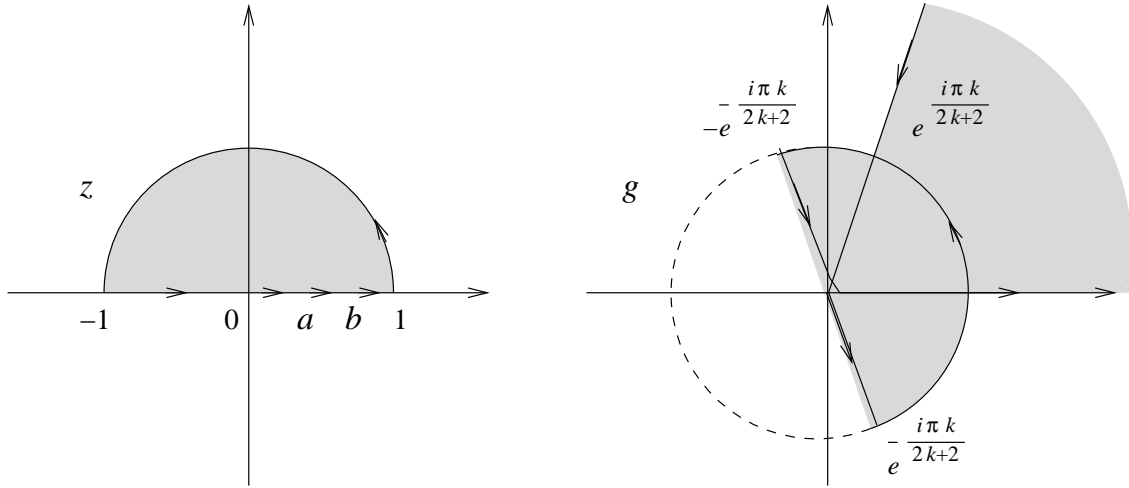}}
\caption{The domain $\Dd$ and its image under $g$.}
\end{figure}

Now consider the following stretches of $\deh\Dd$: $[0,a]$, $[a,b]$, $[b,1]$, $e^{i[0,\pi]}$, $[-1,0]$. After analysing $dh$ and $\tg$ on each of them, we conclude that the projection of $\tX(\deh\Dd)$ on the plane $x_2x_3$, which we call $\Cc$, will be a curve like the ones depicted in Figure 11. In fact, this holds for $k>1$. For $k=1$, the stretch of $\Cc$ that connects $0$ with a point in $Ox_3$ is contained in that axis.  
\\
\input epsf
\begin{figure}[ht]
\centerline{
\epsfxsize 15cm
\epsfbox{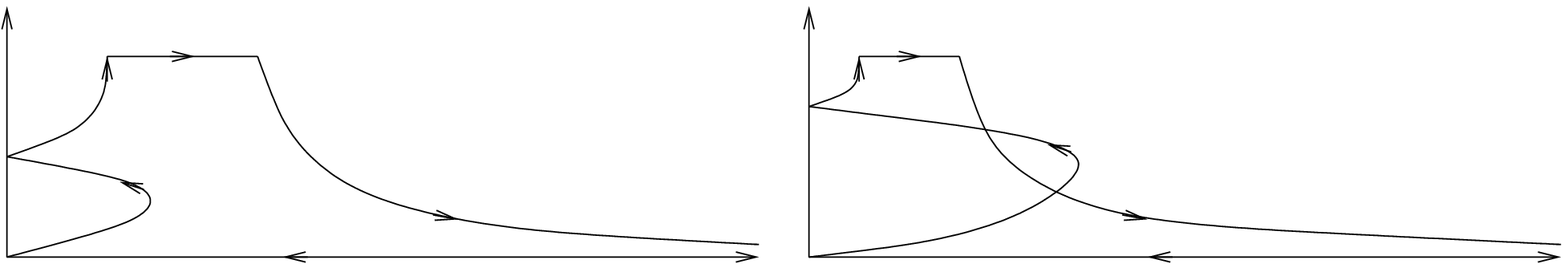}}
\caption{Possible projections of $\tX(\deh\Dd)$ onto $x_2x_3$.}
\end{figure}

Notice that $\tX(\Dd)$ is contained in its convex hull (see [{\bf 22}] for details), which is a subset of $\F:=\{(x_1,x_2,x_3)\in\R^3:x_2\ge 0,0\le x_1\le x_2\tan(\pi/(k+1))$ and $0\le x_3\le 1/2\}$. Successive reflections in the components of $\deh F$ will tessellate $\R^3/\T$ in exactly $8k+8$ congruent pieces. Let $\F'$ be any of these pieces. Of course, $\F\supset\tX(\deh\Dd\setminus\{0\})$, this one all made up by symmetry curves and one ray, in such a way that $\tX(\Ns)\cap\F'\cong\tX(\overline{\Dd}^*)$.  
\\

From (6), the divisor of $g$ is given by
\[
   [g]=\frac{p_1^{k+2}(q_2r_1r_3)^k}{p_2^{k+2}(q_1r_2r_4)^k},
\]
thus
\[
   [dg]=\frac{p_1^{k+1}(q_2r_1r_3)^{k-1}}{p_2^{k+3}(q_1r_2r_4)^{k+1}}\cdot D,
\]
where $D$ is a divisor with $8k+8$ zeroes, because deg$(dg)=-\chi(\overline{\Ns})=8k$. Now, if we had $dg\ne 0$ in $\Dd$, then $g|_\Dd$ would be an unbranched covering. However, any curve in $\Dd$ connecting $0$ and $b$ is taken to a loop with base-point $0$. But $g(\Dd)$ is simply connected, hence the pre-image of any loop should be a closed curve. This contradiction implies that $dg$ has at least one zero in $\Dd$. Now recall that $\R^3/\T$ is tessellated by exactly $8k+8$ pieces congruent to $\F$. Consequently, $dg\ne 0$ on $\deh\Dd\setminus\{0,b\}$, and therefore $\Cc$ is a monotone curve. In particular, all of its stretches depicted in Figure 11 {\it are} convex.
\\

Since $(x_2,x_3):\Dd\to\R^2$ is an immersion and also injective on $\deh\Dd\setminus\{0\}$, then it is a covering map. From [{\bf 10}] or [{\bf 23}], the Gaussian curvature of minimal surfaces is given by
\[
   K=\frac{-16}{(|g|+|g|^{-1})^4}\biggl|\frac{dg/g}{dh}\biggl|^2.
\]
Since neither $g$ nor $dh$ vanishes on $(a,b)$, this means that $K\ne 0$ on this stretch. In Figure 11, this means that $(x_2,x_3)(\Dd)$ contains an open neighbourhood at the right-hand side of the corresponding stretch for $\Cc$. But $(x_2,x_3):\Dd\to\R^2$ is an immersion, which implies that $\Cc$ must be simple. Therefore, Int $\Cc$ is simply connected, and the covering map $(x_2,x_3)|_\Dd$ must be a graph. Consequently, $\tX:\Dd\to\R^3/\T$ is an embedded piece, and from the above discussion of $\F$, $\tX:\Ns\to\R^3/\T$ is an embedding. This concludes the proof of Theorem 1.1.
\\
\\
\centerline{\bf References}
\begin{description}
\itemsep = 0.0 pc
\parsep  = 0.0 pc
\parskip = 0.0 pc
\item{\eh[1]} C.J. Costa, {\it Uniqueness of minimal surfaces embedded in $\R^3$ with total curvature $12\pi$}, J. Differential Geom. {\bf 30} (1989) 597--618.
\item{\eh[2]} M. Callahan, D. Hoffman \& H. Karcher, {\it A family of singly periodic minimal surfaces invariant under a screw motion}, Experiment. Math. {\bf 2} (1993) 157--182.
\item{\eh[3]} M. Callahan, D. Hoffman \& W.H. Meeks, {\it Embedded minimal surfaces with an infinite number of ends}, Invent. Math. {\bf 96} (1989) 459--505.
\item{\eh[4]} M. Callahan, D. Hoffman \& W.H. Meeks, {\it The structure of singly-periodic minimal surfaces}, Invent. Math. {\bf 99} (1990) 455--481.
\item{\eh[5]} H.I. Choi, W.H. Meeks \& B. White, {\it A rigidity theorem for properly embedded minimal surfaces in $\R^3$}, J. Differential Geom. {\bf 32} (1990) 65--76.
\item{\eh[6]} L. Ferrer \& F. Mart\ih n, {\it Minimal surfaces with helicoidal ends}, Math. Z. {\bf 250} (2005) 807--839.
\item{\eh[7]} D. Hoffman \& H. Karcher, {\it Complete embedded minimal surfaces of finite total curvature}, Encyclopedia of Math. Sci., Springer Verlag {\bf 90} (1997) 5--93.
\item{\eh[8]} D. Hoffman, H. Karcher \& F. Wei, {\it The singly periodic genus-one helicoid}, Comment. Math. Helv. {\bf 74} (1999) 248--279.
\item{\eh[9]} D. Hoffman \& W.H. Meeks, {\it Embedded minimal surfaces of finite topology}, Ann. of Math. {\bf 131} (1990) 1--34.
\item{[10]} H. Karcher, {\it Construction of minimal surfaces}, Surveys in Geometry, University of Tokyo (1989) 1--96 and Lecture Notes {\bf 12} (1989), SFB256, Bonn.
\item{[11]} N. Kapouleas, {\it Complete embedded minimal surfaces of finite total curvature}, J. Differential Geom. {\bf 47} (1997) 95--169. 
\item{[12]} L. Lazard-Holly \& W.H. Meeks, {\it Classification of doubly-periodic minimal surfaces of genus zero}, Invent. Math. {\bf 143} (2001) 1--27.
\item{[13]} F.J. L\'opez \& F. Mart\ih n, {\it Complete minimal surfaces in $\R^3$}, Publ. Mat. {\bf 43} (1999) 341--449.
\item{[14]} F.J. L\'opez \& A. Ros, {\it On embedded complete minimal surfaces of genus zero}, J. Differential Geometry {\bf 33} (1991) 293--300.
\item{[15]} F. Mart\ih n, {\it A note on the uniqueness of the periodic Callahan-Hoffman-Meeks surfaces in terms of their symmetries}, Geom. Dedicata {\bf 86} (2001) 185--190.
\item{[16]} F. Mart\ih n \& D. Rodr\ih guez, {\it A characterization of the periodic Callahan-Hoffman-Meeks surfaces in terms of their symmetries}, Duke Math. J. {\bf 89}  (1997) 445--463.
\item{[17]} F. Mart\ih n \& M. Weber, {\it On properly embedded minimal surfaces with three ends}, Duke Math. J. {\bf 107} (2001), 533--559.
\item{[18]} W.S. Massey, {\it Algebraic topology: an introduction}, Graduate Texts in Mathematics, Springer, New York (1967).
\item{[19]} W.H. Meeks, J. P\'erez \& A. Ros, {\it Uniqueness of the Riemann minimal examples}, Invent. Math. {\bf 131} (1998) 107--132.
\item{[20]} W.H. Meeks \& H. Rosenberg, {\it The uniqueness of the helicoid}, Ann. of Math. {\bf 161} (2005) 727--758.
\item{[21]} J.C.C. Nitsche, {\it Lectures on minimal surfaces}, Cambridge University Press, Cambridge (1989).
\item{[22]} R. Osserman, {\it The convex hull property of immersed manifolds}, J. Differential Geometry {\bf 6} (1971/72) 267--270.
\item{[23]} R. Osserman, {\it A survey of minimal surfaces}, Dover, New York, 2nd ed (1986).
\item{[24]} J. P\'erez, M. Rodr\ih guez \& M. Traizet, {\it The classification of doubly periodic minimal tori with parallel ends}, J. Differential Geom. {\bf 69} (2005) 523--577.
\item{[25]} J. P\'erez \& M. Traizet, {\it The classification of singly periodic minimal surfaces with genus zero and Scherk type ends}, Trans. Amer. Math. Soc. (to appear).
\item{[26]} V. Ramos Batista, {\it A family of triply periodic Costa surfaces}, Pacific J. Math. {\bf 212} (2003) 347--370.
\item{[27]} V. Ramos Batista, {\it Noncongruent minimal surfaces with the same symmetries and conformal structure}, Tohoku Math. J. {\bf 56} (2004) 237--254.
\item{[28]} R. Schoen, {\it Uniqueness, symmetry and embeddedness of minimal surfaces}, J. Differential Geom. {\bf 18} (1983) 701--809.
\item{[29]} M. Traizet, {\it An embedded minimal surface with no symmetries}, J. Differential Geom. {\bf 60} (2002) 103--153.
\item{[30]} M. Weber, {\it A Teichm\"uller theoretical construction of high genus singly periodic minimal surfaces invariant under a translation}, Manuscripta Math. {\bf 101} (2000) 125--142.
\item{[31]} M. Wohlgemuth, {\it Minimal surfaces of higher genus with finite total curvature}, Arch. Rational Mech. Anal. {\bf 137} (1997) 1--25.
\end{description}

\hfill IME - U{\fns NIVERSITY OF} S{\fns\~AO} P{\fns AULO}

\hfill R{\fns UA DO} M{\fns AT\~AO} 1010

\hfill 05508-090, S{\fns\~AO} P{\fns AULO} - SP, B{\fns RAZIL}

\hfill {\it E-mail address}: paqs@ime.usp.br
\\

\hfill CMCC - ABC F{\fns EDERAL} U{\fns UNIVERSITY} 

\hfill R{\fns UA} C{\fns ATEQUESE} 242, 3rd F{\fns LOOR} 

\hfill 09090-400, S{\fns ANTO} A{\fns NDR\'E} - SP, B{\fns RAZIL}

\hfill {\it E-mail address}: valerio.batista@ufabc.edu.br
\end{document}